\documentclass[a4paper,11pt,reqno]{amsart}
\usepackage{amsmath}
\usepackage{amsthm}

\usepackage{amssymb}

\usepackage{color}

\usepackage{color}
\usepackage[italian,english]{babel}

\usepackage[utf8x]{inputenc}
\usepackage{fancyhdr}

\usepackage{calc}

\usepackage[margin=2.6cm]{geometry}

\usepackage{srcltx}

\usepackage[T1]{fontenc}

\makeatletter
\def\cleardoublepage{\clearpage\if@twoside \ifodd\c@page\else%
		 \hbox{}%
	 \thispagestyle{empty}
	 \newpage%
	 \if@twocolumn\hbox{}\newpage\fi\fi\fi}
\makeatother

\let\cleardoublepage\clearpage

\addto\captionsitalian{}

\newtheorem{thm}{Theorem}[section]
\newtheorem{cor}[thm]{Corollary}
\newtheorem{lem}[thm]{Lemma}
\newtheorem{pro}[thm]{Proposition}
\newtheorem{den}[thm]{Definition}
\newtheorem{oss}[thm]{Remark}
\numberwithin{equation}{section}

\begin{document}

\title[Porous media equations with homogeneous Neumann boundary conditions]{Sharp short and long time $\mathbf L^{\boldsymbol \infty}$ bounds\\ for solutions
to porous media equations\\ with Neumann boundary conditions}

\author {Gabriele Grillo, Matteo Muratori}

\address {Gabriele Grillo, Matteo Muratori: Dipartimento di Matematica, Politecnico di Milano, Piaz\-za Leonardo da Vinci 32, 20133 Milano, Italy}

\email {gabriele.grillo@polimi.it}

%

\email {matteo1.muratori@mail.polimi.it}

\begin{abstract}
We study a class of nonlinear diffusion equations whose model is the classical porous media equation on domains $\Omega\subseteq{\mathbb R}^N$, $N\ge3$, with homogeneous Neumann boundary conditions. Firstly we improve some known results in such model case, both as concerns sharp $L^{q_0}$-$L^\infty$ regularizing properties of the evolution for short time and as concerns sharp long time asymptotics in the sense of $L^\infty$ convergence of solutions to their mean value. The generality of the discussion allows to consider, almost at the same time,  also weighted versions of the above equation provided an appropriate weighted Sobolev inequality is required to hold. \normalcolor In fact, we show that the validity of a slightly weaker functional inequality  is equivalent to the validity of a suitable $L^{q_0}$-$L^\infty$ bound for solutions to the associated weighted porous media equation.
The long time asymptotic analysis relies as well on the assumed weighted Sobolev inequality only, and allows to prove uniform convergence to the mean value, with the rate predicted by linearization, in such generality. This fact was not known even for the explicit classes of weights previously considered in the literature. 
\end{abstract}

\maketitle

\begin{section}{Introduction}
The aim of this note is to provide some new $L^{q_0}$-$L^\infty$ regularity and asymptotic estimates for solutions to nonlinear diffusion equations whose model is the following \emph{Porous Media Equation} with homogeneous Neumann Boundary conditions:
\begin{equation} \label{eq: pmeNeu}
\begin{cases}
u_t = \Delta\left(u^m\right)  & \textnormal{in} \  \Omega\times(0,\infty)  \\
 \frac{\partial{(u^m)}}{\partial{\mathbf{n}}} =0 & \textnormal{on} \ \partial\Omega \times (0,\infty) \\
 u(\cdot,0)=u_0(\cdot) & \textnormal{in} \ \Omega
 \end{cases} \, ,
\end{equation}
where $\Omega$ is a bounded regular (say $C^1$) domain of $\mathbb{R}^N$, $m>1, N\ge3$ and, as usual, for $y \neq 0$ we define $y^m:=|y|^{m-1} y$ when dealing with non-necessarily positive solutions. Even in this case, in which we shall prove results which will turn out to be sharp both for short and long time, such estimates will improve those obtained in the pioneering paper \cite{AR81} and, for larger classes of data, in  \cite{BG05m} (see \cite{AT, ACLT} for the Neumann problem in the case in which $\Omega$ has infinite measure).  Actually, almost all of these results will also work for the \emph{Weighted Porous Media Equation} introduced in \cite{GMP} (see \cite{KR1, KR2, Eid90, EK, RV08, RV09, KRV10} and references quoted therein for the one-weight case):
\begin{equation} \label{eq: pmeNeuWei}
\begin{cases}
\rho_{\nu} \, u_t =\operatorname{div}\left(\rho_{\mu} \, \nabla{\left( u^m \right)} \right) & \textnormal{in} \  \Omega\times(0,\infty)  \\
 \rho_{\mu} \, \frac{\partial{(u^m)}}{\partial{\mathbf{n}}} =0 & \textnormal{on} \ \partial\Omega \times (0,\infty) \\
 u(\cdot,0)=u_0(\cdot) & \textnormal{in} \ \Omega
 \end{cases} \, .
\end{equation}
Indeed, the only relevant assumption for our method to work is that the weights $\rho_\nu , \rho_\mu $ are strictly positive in $\Omega$ (but may be \it degenerate \rm or \it singular \rm at the boundary), sufficiently regular and such that the Sobolev-type inequality (assuming throughout the paper that $\nu(\Omega)<\infty$)
\begin{equation}\label{eq: soboIntroMedWei}
\left\| v-\overline{v} \right\|_{2\sigma;\nu} \leq C_S \left\| \nabla{v} \right\|_{2;\mu} \ \ \ \forall v \in W^{1,2}(\Omega;\nu,\mu)
\end{equation}
holds true for suitable $C_S>0$, $\sigma > 1$, where $\overline{v}$ is the mean value of $v$ with respect to the measure $\mathrm{d}\nu=\rho_\nu \mathrm{d}\mathbf{x}$ (we also set d$\mu=\rho_\mu {\rm d}\mathbf{x}$) and  $\|\cdot\|_{p;\lambda}$ will denote, for any $p\in[1,\infty)$, the $L^p$ norm w.r.t.\ a given measure $\lambda$ on $\Omega$. In fact the validity of the weaker inequality
\begin{equation}\label{eq: sobDebWei}
\left\|  v \right\|_{2\sigma;\nu} \leq C_S^\prime \left( \left\| \nabla{v}  \right\|_{2;\mu} + \left\| v \right\|_{1;\nu} \right) \ \ \ \forall v \in W^{1,2}(\Omega;\nu,\mu)
\end{equation}
will be sufficient for most of our purposes, and actually will turn out to be \emph{equivalent} to the following bounds for the $L^\infty$ norm of solutions that we shall prove:
\begin{equation*}\label{eq: soboStimaRegWei}
\left\| u(t) \right\|_{\infty} \leq    K \left(t^{-\frac{\sigma}{(\sigma-1)q_0+\sigma(m-1)}} \left\|  u_0 \right\|_{q_0;\nu^{\phantom{a}}}^{\frac{(\sigma-1)q_0}{(\sigma-1)q_0+\sigma(m-1)}} +  \left\| u_0 \right\|_{q_0;\nu}  \right)  \ \ \ \forall t>0 \, ,
\end{equation*}
$K$ being a constant which depends only on $q_0 \in [1,\infty)$, $m$, $C_S$, $\sigma$ and $\nu(\Omega)$.

We shall discuss these issues in greater detail in Section \ref{sec: stimeWei}. One should anyway compare these kind of results with the ones of \cite{GMP}, which involved problems similar to \eqref{eq: pmeNeuWei} under the weaker requirement that the Poincar\'e inequality
\begin{equation}\label{eq: poinIntroMedWei}
\left\| v-\overline{v} \right\|_{2;\nu} \leq C_P \left\| \nabla{v} \right\|_{2;\mu} \ \ \ \forall v \in W^{1,2}(\Omega;\nu,\mu)
\end{equation}
holds true, so that $\min\left[\mathcal{S}\left(-\Delta_{\mathcal N}\right)\setminus\{0\}\right]\ge 1/{C_P^2}>0$, where $\mathcal{S}\left(-\Delta_{\mathcal N}\right)$ denotes the $L^2(\Omega; \nu)$ spectrum of (minus) the weighted Laplacian $-\Delta_{\mathcal N} (v)=-\rho_{\nu}^{-1}\operatorname{div}\left(\rho_{\mu} \, \nabla{ v } \right)$ with homogeneous Neumann boundary conditions. While of course the results of \cite{GMP} continue to hold in the present situation, stronger ones will be shown to be valid and, in particular, solutions corresponding to $L^{q_0}$ data become instantaneously bounded with quantitative bounds on the $L^\infty$ norm, a fact which need not be true under the sole condition \eqref{eq: poinIntroMedWei}, as shown in \cite{GMP}.

As for the long time asymptotics, we discuss first the non-weighted case and show uniform convergence of solutions to the mean value of the initial datum with a sharp rate. Again this bounds improve on the results of \cite{BG05m}. In the weighted case, we shall prove that uniform convergence to the (weighted) mean value of the initial datum occurs, for general $L^1$ data, as a consequence of the sole validity of the Sobolev-type inequality \eqref{eq: soboIntroMedWei}: the proof of this fact is completely different from the one known in the non-weighted context and requires a much more delicate functional argument. In particular, if the weighted mean is not zero, the rate we give is exactly the one predicted by linearization. The existing results (see \cite{KR2, GHP, EiK}) show, in the explicit class of weights considered there, only pointwise or local uniform convergence; however we stress that the Sobolev-type inequality we require need not always hold in their setting. 
As already mentioned, we comment that uniform convergence does not necessarily occur if one assumes only the validity of \eqref{eq: poinIntroMedWei}, as shown in \cite{GMP} by explicit counterexamples. 

For a thorough study of smoothing and decay properties of solutions to large classes of nonlinear evolution equations on ${\mathbb R}^N$ see the monograph \cite{V}, whilst for other work specifically concerned with the connection between functional inequalities and asymptotic properties of solutions to weighted porous media equations we refer the reader to \cite{DG+08, DNS08, W}, remarking however that $L^{q_0}$-$L^\infty$ regularization properties of the evolutions considered are not addressed there. The pioneering papers \cite{KR2}, \cite{GHP} show local uniform convergence of solutions to their mean in some explicit one dimensional, single weight case but do not deal with rates of convergence. For the Neumann problem in the case in which $\Omega$ has \emph{infinite} measure (where convergence to zero is considered) see e.g.\ \cite{AT, ACLT}.

The proofs of existence and uniqueness of weak solutions to problem \eqref{eq: pmeNeu} can be found in \cite{AR81}, \cite{Vaz07}. A proof of existence of solutions (and of uniqueness in the class of {\it energy solutions}) to \eqref{eq: pmeNeuWei} has been given in \cite{GMP}, under certain regularity conditions on the weights: in particular the problem is well-posed provided the weights are sufficiently smooth and locally bounded away from zero in $\Omega$ (while they can be singular or degenerate at the boundary), a fact that we shall assume hereafter.

Our first main result, Theorem \ref{teo: soboTeoReg}, gives a bound on the $L^\infty$ norm of solutions to \eqref{eq: pmeNeu} for short times and $L^{q_0}$ ($q_0 \ge 1$) initial data. This bound improves considerably Theorem 1.1 of \cite{BG05m} (see the discussion in Remark \ref{oss: conf}) and it is sharp in the sense that it captures the exact explosion rates of the well-known Barenblatt solutions \cite{Vaz07} for short time. Bounds for large times are given in Theorem \ref{eq: soboAbsbMNullaPme} for zero-mean data, they being sharp (as already noticed by \cite{AR81} for $L^\infty$ data) and improving the corresponding results of \cite{BG05m}: in particular, an absolute bound for the $L^\infty$ norm follows. The case of data having non-zero mean is also considered in Theorem \ref{eq: teoAsiMNon}, where we show that convergence to the mean value occurs with an exponential rate matching the one predicted by linearization. Once again we improve the results of \cite{BG05m} and in Proposition \ref{lower} we show the sharpness of our bounds by producing a class of data for which a matching lower bound for the rate of convergence holds.

The corresponding results for solutions to \eqref{eq: pmeNeuWei} are given in Section \ref{sec: stimeWei}, see in particular Theorems \ref{teo: soboTeoRegWei}, \ref{cor: soboMNAbsbWei}, \ref{teo: convUnif} and \ref{th: teoAsiMNonPesi}.  Besides proving such results, which are parallel to those obtained for the non-weighted case, we show in Theorem \ref{thm: impInv} and in Corollary \ref{cor: equivNeu} that the validity of a short time $L^{q_0}$-$L^\infty$ regularizing effect for solutions to \eqref{eq: pmeNeuWei} implies in turn the validity of a suitable Sobolev-type inequality like \eqref{eq: sobDebWei}, thus giving a converse to the corresponding regularization theorem and hence showing its optimality, in some sense. Uniform convergence to the (non-zero) weighted mean of the initial datum is shown in Theorem \ref{teo: convUnif}, whereas Theorem \ref{th: teoAsiMNonPesi} gives an explicit rate of convergence.

Finally, Section \ref{sec: exaSobo} provides a concise list of explicit classes of weights for which suitable Sobolev inequalities hold. We remark that almost none of the corresponding nonlinear diffusions seems to have been studied in the existing literature.
\end{section}

\begin{section}{The concept of solution}\label{sec:def}
For the reader's convenience we begin recalling from \cite{GMP} (see also \cite{AR81} and \cite{Vaz07} for the non-weighted case) the concept of solution we are going to consider. We write it in the weighted case only, but clearly the definition applies for the non-weighted case as well just by setting $\rho_\nu, \rho_\mu\equiv1$. Recall that $\mathrm{d}\nu=\rho_\nu \mathrm{d}\mathbf{x}$ and d$\mu=\rho_\mu {\rm d}\mathbf{x}$.

\begin{den}\label{den: solDebNeu}
A function
$$u \in L^2((0,T);L^2(\Omega;\nu)): \ \nabla{(u^m)} \in L^2((0,T);[L^2(\Omega;\mu)]^N) \ \ \ \forall T>0  $$
is a weak solution of \eqref{eq: pmeNeuWei} with initial datum $u_0 \in L^2(\Omega;\nu)$ if it satisfies
\begin{equation*}\label{eq: solDebNeu}
\int_0^T \! \! \int_{\Omega} u(\mathbf{x},t) \eta_t(\mathbf{x},t) \, \mathrm{d}\nu \, \mathrm{d}t = -\int_\Omega u_0(\mathbf{x}) \eta(\mathbf{x},0) \,  \mathrm{d}\nu + \int_0^T \! \! \int_{\Omega} \nabla{\left( u^m\right)}(\mathbf{x},t) \cdot \nabla{\eta}(\mathbf{x},t) \, \mathrm{d}\mu \, \mathrm{d}t
\end{equation*}
$$ \forall \eta \in W^{1,2}((0,T);L^2(\Omega;\nu)): \  \nabla{\eta} \in L^2((0,T);[L^2(\Omega;\mu)]^N)  \, , \ \eta(T)=0 \, . $$
\end{den}
Uniqueness of solutions does not hold in general in the weighted context, even for bounded data and solutions,  however it does indeed  under the additional condition that
\begin{equation*}\label{eq: hpEnergNeu}
u \in L^{m+1}((0,T); L^{m+1}(\Omega;\nu)) \, ,
\end{equation*}
see \cite[Prop. 3.10]{GMP}. We call such a solution, if any, a \emph{weak energy solution} of the equation considered, according to a common terminology. In fact Theorem 3.12 of \cite{GMP} gives a proof of existence of the weak energy solution for initial data $u_0 \in L^{m+1}(\Omega;\nu) $, under some further regularity assumptions on the weights, for example $\rho_\nu, \rho_\mu>0$ and
\begin{equation*}\label{reg weights}
\rho_\nu \in C^{3,\alpha}_{loc}(\Omega) \, , \ \rho_\mu \in C^{2,\alpha}_{loc}(\Omega) \, .
\end{equation*}
For more general data $u_0 \in L^1(\Omega;\nu)$ the correct  extension of the concept of solution is given in \cite[Sect. 3]{GMP} (the discussion being similar to the one first introduced in \cite[Sect. 6.1]{Vaz07} for the non-weighted problem). Such solutions are called \emph{limit solutions}, since they are naturally obtained as limits of energy solutions, that is by approximating $u_0$ with a sequence of data $u_{0n} \in L^{m+1}(\Omega;\nu)$ and exploiting the fundamental $L^1$-comparison principle, namely the fact that if $u$ and $v$ are the solutions to \eqref{eq: pmeNeuWei} with initial data respectively $u_0$ and $v_0$ then
\begin{equation}\label{eq: princConf}
\left\|(u(t)-v(t))_{+} \right\|_{1;\nu} \leq \left\|(u_0-v_0)_{+} \right\|_{1;\nu} \ \ \ \forall t>0 \, .
\end{equation}

Thanks to \eqref{eq: princConf}, in particular, the sequence of energy solutions $\{u_n\}$ corresponding to $\{u_{0n}\}$ is Cauchy in $L^{\infty}((0,\infty);L^1(\Omega;\nu))$. For more details about the mentioned existence results of energy solutions, we refer the reader to \cite[Th. 3.12]{GMP} or to the original \cite[Th. 11.2]{Vaz07} in the non-weighted case.

Throughout the whole paper, when referring to ``the solution'' to the equation at hand, we shall mean without further comment the unique weak energy or limit solution constructed as above. Also, we shall often make use of two fundamental properties of such solutions. In first place, the mean value of the initial datum is preserved along the evolution: see \cite[Prop. 3.13]{GMP} for a proof. In second place, the inequality $\|u(t)\|_{p;\nu} \le \|u(s)\|_{p;\nu}$ holds true for all $p\in[1,\infty]$ and all $t\ge s$, as a consequence of the energy estimates of \cite[Sect. 3]{GMP}; we shall refer to the latter property as ``non-expansivity'' of the norms. For the corresponding results in the non-weighted case see \cite[Th. 11.10]{Vaz07}.
\normalcolor

In the following sections, when dealing with \eqref{eq: pmeNeu}, $\Omega$ will always be a bounded $C^1$ domain of $\mathbb{R}^N$. As for the weighted case \eqref{eq: pmeNeuWei} regularity of $\Omega$ in principle is not needed, although it may be hidden in the assumed validity of the corresponding weighted Sobolev inequalities.
\end{section}

\begin{section}{Short time $L^{q_0}$-$L^\infty$ bounds}\label{sec: stimeReg}
We start recalling a useful numerical lemma given in \cite{GMP}.
\begin{lem} \label{eq: lemminoNumerico}
Given $\alpha, \beta \in (0,1)$, with $\alpha > \beta$, there exists a constant $c=c(\alpha, \beta)>0$ such that
\begin{equation*} \label{eq: minoreLemmino}
x^{-\alpha}y^{1-\alpha}+x^{-\beta}y^{1-\beta}+y \leq c(\alpha,\beta)(x^{-\alpha}y^{1-\alpha}+y) \ \ \ \forall x,y \in \mathbb{R}^{+}  \, .
\end{equation*}
\begin{proof}
We need to show that
\begin{equation*} \label{eq: R}
R(x,y)=\frac{x^{-\beta}y^{1-\beta}} {x^{-\alpha}y^{1-\alpha}+y }
\end{equation*}
is bounded in $\mathbb{R}^{+} \! \times \mathbb{R}^{+}$ by a constant which depends only on $\alpha$ and $\beta$, and in order to do that one finds explicitly the zeros of $R_{x}(\cdot,y)$ for any given $y$.
\end{proof}
\end{lem}

Now we prove one of the main results of this paper which, as we shall remark later, is a sharp improvement of the $L^{q_0}$-$L^\infty$ regularity estimate first provided by Theorem 1.1 of \cite{BG05m}.

\begin{thm}\label{teo: soboTeoReg}
Let $u$ be the solution of \eqref{eq: pmeNeu} corresponding to an initial datum $u_0 \in L^{q_0}(\Omega)$, with $q_0 \in [1,\infty)$. The following estimate holds:
\begin{equation}\label{eq: soboStimaReg}
\left\| u(t) \right\|_{\infty} \leq    K \left(t^{-\frac{N}{2q_0+N(m-1)}} \left\|  u_0 \right\|_{q_0^{\phantom{a}}}^{\frac{2q_0}{2q_0+N(m-1)}} +  \left\| u_0 \right\|_{q_0^{\phantom{a}}}  \right)  \ \ \ \forall t>0 \, ,
\end{equation}
where $K$ is a constant which depends on $m$, $C_S$, $N$, $|\Omega|$ and can be taken to be independent of $q_0$.
\begin{proof}
We shall proceed by means of a classical Moser iterative technique. Firstly we consider an initial datum $u_0 \in L^\infty(\Omega)$: the fact that the estimate we shall obtain will not depend on $\|u_0\|_\infty$ will allow us to extend it to general $L^{q_0}$ data thanks to a well-known argument that we shall recall at the end of this proof.

Given $t>0$, let us consider the sequence of time steps $t_k=t \, (1-2^{-k}) $. Clearly, $t_0=0$ and $t_{\infty}=t$. Also, let $\{p_k\}$ be an increasing sequence of positive numbers such that $p_0=q_0$ and $p_\infty=\infty$, which we shall explicitly define later. For the moment, we assume in addition that $q_0 \in (1,\infty) \cap [m-1,\infty)$ (afterwards we shall be able to remove this hypothesis). Multiplying \eqref{eq: pmeNeu} by $u^{p_k - 1}$ and integrating in $\Omega \times (t_k,t_{k+1})$ (and neglecting $\| u(t_{k+1}) \|_{p_k}^{p_k}$) gives
\begin{equation} \label{eq: soboRegDecre1}
\frac{4 (p_k-1)p_k m}{(p_k+m-1)^2} \int_{t_k}^{t_{k+1}} \! \! \int_\Omega \left|\nabla{\left(u^{\frac{p_k+m-1}{2}} \right)}(\mathbf{x},s)\right|^2 \,   \mathrm{d}\mathbf{x}  \, \mathrm{d}s \leq  \left\| u(t_k) \right\|_{p_k}^{p_k}  \, .
\end{equation}
With no loss of generality, suppose $|\Omega|=1$. In order to suitably handle the left hand side of \eqref{eq: soboRegDecre1}, it is convenient to notice that the validity of the Sobolev inequality
\begin{equation*}\label{eq: soboIntroMed}
\left\| v-\overline{v} \right\|_{\frac{2N}{N-2}} \leq C_S \left\| \nabla{v} \right\|_{2} \ \ \ \forall v \in W^{1,2}(\Omega)
\end{equation*}
implies the validity of the inequality
\begin{equation}\label{eq: soboRisc1}
\frac{1}{2 C_S^2} \left\| v \right\|_{2\sigma}^2 - \frac{1}{C_S^2} \left\| v \right\|_{1}^2 \leq \left\| \nabla{v} \right\|_{2}^2 \ \ \ \forall v \in W^{1,2}(\Omega) \, ,
\end{equation}
where $\sigma={N}/(N-2)$. Upon applying \eqref{eq: soboRisc1} to the function $u^{(p_k+m-1)/{2}}$ in \eqref{eq: soboRegDecre1}, we get:
\begin{equation*} \label{eq: soboRegDecre2}
\begin{split}
& \frac{2(p_k-1)p_k m}{C_S^2(p_k+m-1)^2} \int_{t_k}^{t_{k+1}} \left\| u(s) \right\|_{\sigma(p_k+m-1)}^{p_k+m-1} \, \mathrm{d}s  \leq \\
& \leq \left\| u(t_{k}) \right\|_{p_k}^{p_k}  +  \frac{4(p_k-1)p_k m}{C_S^2 (p_k+m-1)^2} \int_{t_k}^{t_{k+1}} \left\| u(s) \right\|_{\frac{p_k+m-1}{2}}^{p_k+m-1} \, \mathrm{d}s   \, .
\end{split}
\end{equation*}
Since $q_0 \geq m-1$ and $p_k$ is increasing, we can control $\| u \|_{(p_k+m-1)/{2}}$ with $\| u \|_{p_k}$. By that and by the non-expansivity of the norms, we deduce:
\begin{equation} \label{eq: soboRegDecre3}
\frac{(p_k-1)p_k m}{C_S^2(p_k+m-1)^2} \, 2^{-k} \, t \left\| u(t_{k+1}) \right\|_{p_{k+1}^{\phantom{a}}}^{\frac{p_{k+1}}{\sigma}} \!\! \leq  \left\| u(t_{k}) \right\|_{p_k}^{p_k}  +  \frac{2 (p_k-1)p_k m}{C_S^2 (p_k+m-1)^2} \, 2^{-k} \, t \left\| u(t_{k}) \right\|_{p_k}^{p_k+m-1}   \, ,
\end{equation}
provided $p_{k+1}=\sigma(p_k+m-1)$. Now we assume $\| u_0 \|_\infty=1$. This hypothesis, again together with the non-expansivity of the norms, ensures that
$$\| u(t_k) \|_{p_k}^{p_k+m-1} \leq \| u(t_k) \|_{p_k}^{p_k} \, , $$
so that \eqref{eq: soboRegDecre3} reads
\begin{equation} \label{eq: soboRegDecre4}
\left\| u(t_{k+1}) \right\|_{p_{k+1}^{\phantom{a}}}^{\frac{p_{k+1}}{\sigma}} \!\! \leq \frac{C_S^2(p_k+m-1)^2}{(p_k-1)p_k m}  \, 2^{k} \, t^{-1} \left\| u(t_{k}) \right\|_{p_k}^{p_k}  +  2 \left\| u(t_{k}) \right\|_{p_k}^{p_k}   \, .
\end{equation}
Clearly there exists a suitable constant $D=D(q_0,m,C_S,N)$ such that \eqref{eq: soboRegDecre4} simplifies to
\begin{equation*} \label{eq: soboSucc1}
\left\| u(t_{k+1}) \right\|_{p_{k+1}} \leq D^{\frac{k+1}{p_{k+1}}}(t^{-1}+1)^{\frac{\sigma}{p_{k+1}}} \left\| u(t_{k}) \right\|_{p_k^{\phantom{a}}}^{\sigma\frac{p_k}{p_{k+1}}}    \, .
\end{equation*}
Setting $U_{k}=\|u(t_{k}) \|_{p_{k}}$, it is straightforward to check that the sequence $\{U_k\}$ satisfies
\begin{equation} \label{eq: soboSucc2}
U_{k+1} \leq D^{\frac{\sigma^{k+2}-(k+2)\sigma+k+1 }{p_{k+1}\left(\sigma-1\right)^2}} (t^{-1}+1)^{\frac{\sigma^{k+2}-\sigma}{p_{k+1}(\sigma-1)} } {U_0}^{q_0 \frac{\sigma^{k+1}}{p_{k+1}}}   \, .
\end{equation}
Also, one can verify that $p_{k}=(q_0-A)\sigma^k+A$, with $A=\frac{\sigma}{\sigma-1}(1-m)=\frac{N}{2}(1-m)$. Letting $k\rightarrow \infty $ in \eqref{eq: soboSucc2}, we get (from now on $D=D(q_0,m,C_S,N)$ will denote a generic constant which may differ from line to line):
\begin{equation} \label{eq: soboStimaNoReg1}
\left\| u(t) \right\|_{\infty}=\lim_{k \rightarrow \infty}{\left\| u(t) \right\|_{p_{k+1}}} \leq  \liminf_{k \rightarrow \infty}{U_{k+1}} \leq  D (t^{-1}+1)^{\frac{N}{2q_0+N(m-1)}} {\left\| u_0 \right\|}_{q_0^{\phantom{a}}}^{\frac{2q_0}{2q_0+N(m-1)}} \, .
\end{equation}
Actually, \eqref{eq: soboStimaNoReg1} is \emph{not} an $L^{q_0}$-$L^\infty$ regularity estimate. Indeed, recall that it has been obtained for initial data $u_0$ such that $\|u_0 \|_{\infty}=1$. Through a simple \emph{time scaling} argument we can deduce from it an estimate for general $L^\infty$ data. That is, given a solution $u(\cdot,t)$ to \eqref{eq: pmeNeu} corresponding to the initial datum $u_0$, it is straightforward to check that $\widehat{u}(\cdot,t)=\frac{1}{\lambda}u(\cdot,\lambda^{1-m} \, t)$ is the solution of the same equation corresponding to the initial datum $\widehat{u}_0=\frac{1}{\lambda}u_0$. Choosing $\lambda=\| u_0 \|_{\infty}$ and applying \eqref{eq: soboStimaNoReg1} to $\widehat{u}(t)$ we conclude that
\begin{equation} \label{eq: soboStimaNoReg2}
\left\| u(t) \right\|_{\infty} = \left\|u_0 \right\|_{\infty} \widehat{u}(\left\| u_0 \right\|_{\infty}^{m-1}\,t) \leq  D (t^{-1}+ \left\| u_0 \right\|_{\infty}^{m-1})^{\frac{N}{2q_0+N(m-1)} } {\left\| u_0 \right\|}_{q_0^{\phantom{a}}}^{\frac{2q_0}{2q_0+N(m-1)}}   \, .
\end{equation}
Estimate \eqref{eq: soboStimaNoReg2}, as such, is not of great interest. However, it is possible to reduce its dependence on $\| u_0 \|_{\infty}$ in the following way. First, let us rewrite it as
\begin{equation*} \label{eq: soboStimaNoReg3}
\left\| u(t) \right\|_{\infty} \leq  D \left(t^{-\frac{N}{2q_0+N(m-1)}} {\left\| u_0 \right\|}_{q_0^{\phantom{a}}}^{\frac{2q_0}{2q_0+N(m-1)}}  + \left\| u_0 \right \|_{\infty^{\phantom{a}}}^{\frac{N(m-1)}{2q_0+N(m-1)}} {\left\| u_0 \right\|}_{q_0^{\phantom{a}}}^{\frac{2q_0}{2q_0+N(m-1)}}  \right)  \, ,
\end{equation*}
that is, setting $\theta=\frac{N(m-1)}{2q_0+N(m-1)}$,
\begin{equation} \label{eq: soboStimaNoReg4}
\left\| u(t) \right\|_{\infty} \leq  D \left( t^{-\frac{\theta}{m-1}} {\left\| u_0 \right\|}_{q_0}^{1-\theta}  + \left\| u_0 \right\|_{\infty}^{\theta} {\left\| u_0 \right\|}_{q_0}^{1-\theta}     \right)   \, .
\end{equation}
Now it is convenient to exploit the semigroup property, by shifting the time origin from $0$ to $t/2$, and the non-expansivity of the $L^{q_0}$ norm. This leads us to
\begin{equation*} \label{eq: soboStimaNoReg5}
\left\| u(t) \right\|_{\infty} \leq  D \left( t^{-\frac{\theta}{m-1}} {\left\| u_0 \right\|}_{q_0}^{1-\theta}  + \left\| u(t/2) \right\|_{\infty}^{\theta} \left\| u_0 \right\|_{q_0}^{1-\theta}     \right)   \, ;
\end{equation*}
applying \eqref{eq: soboStimaNoReg4} to $\| u(t/2) \|_{\infty}$ we obtain
\begin{equation*} \label{eq: soboStimaNoReg6}
\left\| u(t) \right\|_{\infty} \leq  D \left( t^{-\frac{\theta}{m-1}} {\left\| u_0 \right\|}_{q_0}^{1-\theta}  + t^{-\frac{\theta^2}{m-1}} {\left\| u_0 \right\|}_{q_0}^{1-\theta^2} +\left\| u_0 \right\|_{\infty}^{\theta^2} {\left\| u_0 \right\|}_{q_0}^{1-\theta^2}     \right)   \, ;
\end{equation*}
it is then clear that proceeding in this way along $n$ steps one arrives at
\begin{equation} \label{eq: soboStimaNoRegK}
\left\| u(t) \right\|_{\infty} \leq  D(n,\cdot) \left( t^{-\frac{\theta}{m-1}} {\left\| u_0 \right\|}_{q_0}^{1-\theta}  +\ldots+ t^{-\frac{\theta^n}{m-1}} {\left\| u_0 \right\|}_{q_0}^{1-\theta^n} +\left\| u_0 \right\|_{\infty}^{\theta^n} \left\| u_0 \right\|_{q_0}^{1-\theta^n}     \right)   \, .
\end{equation}
In order to remove the dependence on the $L^\infty$ norm, we need a suitable $L^{q_0}$-$L^\infty$ regularity estimate. To this end, one can reason in the following way. Suppose $\| u_0 \|_{q_0} = 1 $. From \eqref{eq: soboRegDecre3}, setting again $U_k=\| u(t_k) \|_{p_k}$, we have:
\begin{equation} \label{eq: soboStimaRegRoz1}
{U_{k+1}}^{\frac{p_{k+1}}{\sigma}}  \leq \frac{C_S^2(p_k+m-1)^2}{(p_k-1)p_k m}  \, 2^{k} \, t^{-1} {U_k}^{p_k}  +  2 {U_k}^{p_k+m-1}   \, .
\end{equation}
Consider a \emph{solution} to \eqref{eq: soboStimaRegRoz1}, namely a sequence $\{V_k\}$ such that $V_0=U_0=1$ and
\begin{equation*} \label{eq: soboStimaRegRoz2}
{V_{k+1}}^{\frac{p_{k+1}}{\sigma}}  = \frac{C_S^2(p_k+m-1)^2}{(p_k-1)p_k m}  \, 2^{k} \, t^{-1} {V_k}^{p_k}  +  2 {V_k}^{p_k+m-1}   \, ;
\end{equation*}
it is easy to verify (by induction) that $U_k \leq V_k $ and $V_k \geq 1$. Proceeding similarly to the case $\| u_0 \|_\infty=1$ (now the leading term on the right hand side is ${V_k}^{p_k+m-1}$) we deduce that for a suitable constant $D=D(q_0,m,C_S,N)$ the sequence $\{V_k\}$ also satisfies the following recursive inequality:
\begin{equation} \label{eq: soboStimaRegRoz3}
{V_{k+1}} \leq D^{\frac{k+1}{p_{k+1}}} \, ( t^{-1} + 1)^{\frac{\sigma}{p_{k+1}}} {V_k}  \, .
\end{equation}
Solving \eqref{eq: soboStimaRegRoz3} and exploiting the fact that $p_k \geq q_0 \, \sigma^k$ one gets (recall that $D$ may vary from line to line)
\begin{equation*} \label{eq: soboStimaRegRoz4}
{V_{k+1}} \leq D^{\frac{(k+2)(1-\sigma)\sigma^{-k-1}-\sigma^{-k-1}+\sigma } {q_0 (\sigma-1)^2} } \, ( t^{-1} + 1)^{\frac{\sigma-\sigma^{-k}}{q_0(\sigma-1) }} \, ,
\end{equation*}
and passing to the limit as $k\rightarrow \infty$,
\begin{equation*} \label{eq: soboStimaRegRoz5}
\left\| u(t) \right\|_{\infty} = \lim_{k \rightarrow \infty}{\left\| u(t) \right\|_{p_k}} \leq \liminf_{k \rightarrow \infty}{U_k} \leq \liminf_{k \rightarrow \infty}{V_k} \leq D \, ( t^{-1} + 1)^{\frac{N}{2 q_0}}  \, .
\end{equation*}
Again, by means of a time scaling argument (now with $\lambda=\| u_0 \|_{q_0}$), we have that the final $L^{q_0}$-$L^\infty$ regularity estimate is
\begin{equation*} \label{eq: soboStimaRegRozFin0}
\left\| u(t) \right\|_{\infty} \leq D \, ( t^{-1} +  \left\| u_0 \right\|_{q_0}^{m-1} )^{\frac{N}{2 q_0}} \left\| u_0 \right\|_{q_0^{\phantom{a}}}^{1-\frac{N(m-1)}{2q_0}}  \, ,
\end{equation*}
which reads, upon setting $\delta={N(m-1)}/{2q_0}$,
\begin{equation} \label{eq: soboStimaRegRozFin1}
\left\| u(t) \right\|_{\infty} \leq D \, \left( t^{-\frac{\delta}{m-1}} \left\| u_0 \right\|_{q_0}^{1-\delta} +  \left\| u_0 \right\|_{q_0} \right) \, .
\end{equation}
Of course the exponents involved in \eqref{eq: soboStimaRegRozFin1} are not very satisfactory, in particular when $q_0$ is small with respect to $N$. However, combining this estimate together with \eqref{eq: soboStimaNoRegK} one can obtain a much better $L^{q_0}$-$L^\infty$ regularity result. Indeed, by means of the usual shift of the time origin to $t/2$ in \eqref{eq: soboStimaNoRegK} and thanks to \eqref{eq: soboStimaRegRozFin1} evaluated at time $t/2$, we have:
\begin{equation*} \label{eq: soboStimaQuasiReg}
\left\| u(t) \right\|_{\infty} \leq  D(\cdot,n) \left( t^{-\frac{\theta}{m-1}} {\left\| u_0 \right\|}_{q_0}^{1-\theta} +\ldots+ t^{-\frac{\theta^n}{m-1}} {\left\| u_0 \right\|}_{q_0}^{1-\theta^n} + t^{-\frac{\delta\theta^n}{m-1}} {\left\| u_0 \right\|}_{q_0}^{1-\delta\theta^n}  +  {\left\| u_0 \right\|}_{q_0} \right)  \, ;
\end{equation*}
now we pick $n$ great enough so that $ \delta \theta^n < \theta$, and apply iteratively Lemma \ref{eq: lemminoNumerico} with $x=t^{1/(m-1)}$, $y=\| u_0\|_{q_0}$, $\alpha=\theta$ and $\beta=\delta\theta^n$ at the first step and then $\beta=\theta^j$ along $j=n\ldots 2$. We thus get
\begin{equation}\label{eq: soboStimaRegSottoQuasiLast}
\left\| u(t) \right\|_{\infty}  \leq D  \left(  t^{-\frac{N}{2q_0+N(m-1)}} \left\| u_0 \right\|_{q_0^{\phantom{a}}}^{\frac{2q_0}{2q_0+N(m-1)}} + \left\| u_0 \right\|_{q_0} \right) \, ,
\end{equation}
which is valid, as stated at the beginning of the proof, for $q_0 \in (1,\infty) \cap [m-1,\infty) $.

We now extend the above estimates to the case of a general $q_0 \geq 1$. In fact, for notational simplicity, we shall give the complete proof only for $q_0=1$.
We shall exploit a technique similar to the one used in the proof of Corollary 8.1 from \cite{PQRV11}. That is, first of all consider the analogue of estimate \eqref{eq: soboStimaRegSottoQuasiLast} in the time interval $[t/2,t]$, for $q_0=m$:
\begin{equation*}\label{eq: soboStimaRegCaso1}
\left\| u(t) \right\|_{\infty}  \leq D  \left(  \left( \frac{t}{2} \right)^{-\frac{N}{2m+N(m-1)}} \left\| u({t}/{2}) \right\|_{m^{\phantom{a}}}^{\frac{2m}{2m+N(m-1)}} + \left\| u({t}/{2}) \right\|_{m} \right) \, ;
\end{equation*}
using the inequality $\| \cdot \|_{m} \le \| \cdot \|_{\infty^{\phantom{a}}}^{(m-1)/m} \| \cdot \|_1^{1/m} $ and the non-expansivity of the $L^1$ norm we obtain
\begin{equation}\label{eq: soboStimaRegCaso2}
\left\| u(t) \right\|_{\infty}  \leq D  \left(  2^{\frac{\gamma}{m-1}} t^{-\frac{\gamma}{m-1}} \left\| u({t}/{2}) \right\|_{\infty^{\phantom{a}}}^{\frac{(m-1)(1-\gamma)}{m}} \left\| u_0 \right\|_{1^{\phantom{a}}}^{\frac{1-\gamma}{m}}  + \left\| u({t}/{2}) \right\|_{\infty^{\phantom{a}}}^{\frac{m-1}{m}} \left\| u_0  \right\|_1^{\frac{1}{m}}  \right) \, ,
\end{equation}
where
$$ \gamma=\frac{N(m-1)}{2m+N(m-1)}  \, .$$
In order to handle \eqref{eq: soboStimaRegCaso2}, we can argue as in the previous part of this proof. That is, consider first an initial datum $u_0$ such that $\| u_0 \|_\infty \leq 1$. This, in particular, implies that $\| u(t/2) \|_\infty \leq 1 $ and $\| u_0 \|_1 \leq 1 $, so that \eqref{eq: soboStimaRegCaso2} becomes 
\begin{equation}\label{eq: soboStimaRegCaso3}
\left\| u(t) \right\|_{\infty}  \leq D \, \left( t^{-\frac{\gamma}{m-1}} + 1 \right) \left\| u({t}/{2}) \right\|_{\infty^{\phantom{a}}}^{\frac{(m-1)(1-\gamma)}{m}} \left\| u_0 \right\|_{1^{\phantom{a}}}^{\frac{1-\gamma}{m}}  \, ,
\end{equation}
up to absorbing $2^{{\gamma}/(m-1)}$ into $D$. It is apparent that iterating estimate \eqref{eq: soboStimaRegCaso3} along $k$ steps one gets
\begin{equation}\label{eq: soboStimaRegCaso4}
\begin{split}
 \left\| u(t) \right\|_{\infty}  \leq &
 2^{\frac{\gamma}{m-1} \sum_{h=0}^{h=k-1}  h \left(\frac{(m-1)(1-\gamma)}{m} \right)^h}  \, \left[D \left( t^{-\frac{\gamma}{m-1}} + 1 \right) \right]^{\sum_{h=0}^{h=k-1} \left(\frac{(m-1)(1-\gamma)}{m} \right)^h } \times \\
 & \times \left\| u\left({t}/{2^k} \right) \right\|_{\infty^{\phantom{a}}}^{\left(\frac{(m-1)(1-\gamma)}{m}\right)^k} \left\| u_0 \right\|_{1^{\phantom{a}}}^{\frac{1-\gamma}{m}\sum_{h=0}^{h=k-1} \left(\frac{(m-1)(1-\gamma)}{m} \right)^h }  \, .
\end{split}
\end{equation}
Passing to the limit in \eqref{eq: soboStimaRegCaso4} as $k \to \infty $ we end up with
\begin{equation}\label{eq: soboStimaRegCaso5}
\left\| u(t) \right\|_\infty \leq D \left( t^{-\frac{1}{m-1} \frac{m\gamma}{1+(m-1)\gamma}} \left\| u_0 \right\|_1^{ \frac{1-\gamma}{1+(m-1)\gamma} }  + \left\| u_0 \right\|_1^{\frac{1-\gamma}{1+(m-1)\gamma} } \right)
\end{equation}
for another suitable constant $D$. Estimate \eqref{eq: soboStimaRegCaso5} is analogous to \eqref{eq: soboStimaNoReg1}. By reasoning likewise we can obtain again \eqref{eq: soboStimaNoRegK}, now with $q_0=1$ and $\theta=\frac{m\gamma}{1+(m-1)\gamma}$:
\begin{equation} \label{eq: soboStimaNoRegK2}
\left\| u(t) \right\|_{\infty} \leq  D(n,\cdot) \left( t^{-\frac{\theta}{m-1}} {\left\| u_0 \right\|}_{1}^{1-\theta}  +\ldots+ t^{-\frac{\theta^n}{m-1}} {\left\| u_0 \right\|}_{1}^{1-\theta^n} +\left\| u_0 \right\|_{\infty}^{\theta^n} \left\| u_0 \right\|_{1}^{1-\theta^n}     \right)   \, .
\end{equation}
In order to remove the dependence of the right hand side of \eqref{eq: soboStimaNoRegK2} on $\| u_0 \|_\infty$, it is convenient to suppose first that $\| u_0 \|_1 \geq 1 $ and look for a suitable $L^{1}$-$L^\infty$ regularity estimate. Indeed, considering again \eqref{eq: soboStimaRegCaso2}, it is easy to check that the worst possible case occurs when $\| u(s) \|_\infty \geq 1 $ for all $s<t$: assuming that, clearly $\| u({t}/{2}) \|_{\infty^{\phantom{a}}}^{(m-1)(1-\gamma)/{m}} \leq  \| u({t}/{2}) \|_{\infty^{\phantom{a}}}^{(m-1)/{m}} $, and since also $\| u_0 \|_{1}^{(1-\gamma)/{m}} \leq  \| u_0 \|_{1}^{1/{m}}$ we have that \eqref{eq: soboStimaRegCaso2} reads
 \begin{equation}\label{eq: soboStimaRegCaso6}
\left\| u(t) \right\|_{\infty}  \leq D \, \left( t^{-\frac{\gamma}{m-1}} + 1 \right) \left\| u({t}/{2}) \right\|_{\infty^{\phantom{a}}}^{\frac{m-1}{m}}  \left\| u_0 \right\|_{1^{\phantom{a}}}^{\frac{1}{m}}  \, .
\end{equation}
From \eqref{eq: soboStimaRegCaso6}, through analogous computations as above we easily get
\begin{equation}\label{eq: soboStimaRegCaso7}
\left\| u(t) \right\|_\infty \leq D \left( t^{-\frac{m\gamma}{m-1}}  \left\| u_0 \right\|_1 + \left\| u_0 \right\|_1 \right) \, ,
\end{equation}
which holds provided $\| u_0 \|_1 \ge 1 $. By means of a time scaling argument (recall how we obtained \eqref{eq: soboStimaRegRozFin1}) \eqref{eq: soboStimaRegCaso7} becomes
\begin{equation}\label{eq: soboStimaRegCaso8}
\left\| u(t) \right\|_\infty \leq D \left( t^{-\frac{m\gamma}{m-1}}  \left\| u_0 \right\|_1^{1-m\gamma} + \left\| u_0 \right\|_1 \right) \, ,
\end{equation}
which is valid for all $u_0 \in L^1(\Omega)$. By choosing $n$ in \eqref{eq: soboStimaNoRegK2} great enough so that $m \gamma \theta^n < \theta $ and then combining \eqref{eq: soboStimaNoRegK2} with \eqref{eq: soboStimaRegCaso8} through the usual $t/2$-shift argument, we obtain:
\begin{equation}\label{eq: soboStimaRegCasoFinale}
\left\| u(t) \right\|_\infty \leq D \left( t^{-\frac{1}{m-1} \frac{m\gamma}{1+(m-1)\gamma}} \left\| u_0 \right\|_1^{ \frac{1-\gamma}{1+(m-1)\gamma} }  + \left\| u_0 \right\|_1\right) \, ,
\end{equation}
which is exactly \eqref{eq: soboStimaRegSottoQuasiLast} with $q_0$ replaced by $1$. As previously stated, a similar strategy also works when $q_0 \in (1,m-1]$, provided this latter interval is not empty. Following the above proof it is easy to realize that the constant $D$ in the inequality corresponding to \eqref{eq: soboStimaRegCasoFinale} when $\|\cdot\|_{q_0}$ instead of $\|\cdot\|_1$ is considered, is bounded as a function of $q_0 \in [1,m-1]$. Moreover, one can notice that the constant $D$ in  \eqref{eq: soboStimaRegSottoQuasiLast} is locally bounded as a function of $q_0 \in (1,\infty) \cap [m-1,\infty) $ and it remains bounded as $q_0\to \infty$. Hence \eqref{eq: soboStimaReg} holds for all $q_0 \ge 1$, with a multiplicative constant $K$ independent of $q_0$.

Finally, we are left with removing the hypotheses $|\Omega|=1$ and $u_0 \in L^\infty(\Omega)$. For the first one it is enough to proceed through a \emph{spatial scaling} argument exactly as explained in the end of the proof of Theorem 5.4 of \cite{GMP}. Namely, if $u(\mathbf{x},t)$ is a solution of \eqref{eq: pmeNeu} in a domain $\Omega$ of measure $|\Omega|$ with initial datum $u_0(\mathbf{x})$ then
\begin{equation}\label{eq: scalingSp}
\widetilde{u}(\widetilde{\mathbf{x}},t)=|\Omega|^{-\frac{2}{N(m-1)}}u\left(|\Omega|^{\frac{1}{N}} \widetilde{\mathbf{x}},t\right)
\end{equation}
is also a solution of \eqref{eq: pmeNeu} in the domain $\widetilde{\Omega}=\Omega/|\Omega|^{\frac{1}{N}}$ of measure $1$ with initial datum $$\widetilde{u}_0(\widetilde{\mathbf{x}})=|\Omega|^{-\frac{2}{N(m-1)}}u_0\left(|\Omega|^{\frac{1}{N}} \widetilde{\mathbf{x}}\right) \, . $$
Noticed that, one applies \eqref{eq: soboStimaReg} to $\widetilde{u}$ and goes back to the original solution $u$ by means of \eqref{eq: scalingSp} and the relations
$$ \left\| \widetilde{u} \right\|_{p}= |\Omega|^{-\frac{2}{N(m-1)}-\frac{1}{p}} \left\| u \right\|_{p} \, , \ C_S(\widetilde{\Omega})=|\Omega|^{\frac{1}{2}-\frac{1}{2\sigma}-\frac{1}{N}}C_S(\Omega) \, . $$
In consequence of that, it all amounts to admit the dependence of the multiplicative constant appearing in \eqref{eq: soboStimaReg} on $|\Omega|$ too.

The extension to general $L^{q_0}$ data is handled by means of a standard argument. That is, given $u_0 \in L^{q_0}(\Omega)$, consider a sequence $\{u_{0n}\} \subset L^\infty(\Omega)$ which converges to $u_0$ in $L^{q_0}(\Omega)$ and the corresponding solutions to \eqref{eq: pmeNeu} $\{u_n\}$ and $u$. From the $L^1$-comparison principle we deduce that, for any given $t>0$, $u_n(t)\rightarrow u(t)$ in $L^1(\Omega)$; moreover from \eqref{eq: soboStimaReg} we also learn that, up to subsequences, $\{u_n(t)\}$ converges in the weak$^\ast$ topology of $L^\infty(\Omega)$ to an element $w \in L^\infty(\Omega)$. The identification between $u(t)$ and $w$ is straightforward, and estimate \eqref{eq: soboStimaReg} is preserved to the limit thanks to the weak$^\ast$ lower semicontinuity of the $L^\infty$ norm.
\end{proof}
\end{thm}
\begin{oss}\label{oss: conf}\rm
We stress that estimate \eqref{eq: soboStimaReg} in fact improves, with respect to the dependence on the time variable and on the parameters $q_0$, $m$, $N$, the one first provided by Theorem 1.1 of \cite{BG05m} (obtained by means of a Gross differential method), which reads
\begin{equation}\label{eq: stimaBG05}
\left\|  u(t) \right\|_\infty \leq C \, t^{-\frac{\alpha}{m-1}} \left\| u_0 \right\|_{q_0}^{1-\alpha} \, e^{E_0 \left\| u_0 \right\|_{1\vee(m-1)}^{m-1} t} \, ,
\end{equation}
where $C,E_0$ are suitable constants that depend only on $q_0$, $m$, $C_S$, $N$, $|\Omega|$ and
\begin{equation}\label{eq: stimaBG05alph}
\alpha=\left[ 1-\left( \frac{q_0}{q_0+m-1} \right)^{\frac{N}{2}} \right] \, .
\end{equation}
Indeed, it is plain that
$$ \alpha > \frac{N(m-1)}{2q_0 + N(m-1) } \ \ \ \forall m>1 \, , \ \forall N > 2 \, , \ \forall q_0 \geq 1 \, .$$
Also, note that \eqref{eq: stimaBG05} blows up as $t \rightarrow \infty$. Moreover, it only holds for initial data belonging to $L^{q_0}(\Omega)$, with $q_0 \geq 1 \vee (m-1) $. Instead, estimate \eqref{eq: soboStimaReg} does not blow up and holds for all initial data which lie in $L^{q_0}(\Omega)$ for any $q_0\ge1$.
\end{oss}

\begin{oss}\label{oss: barenblatt} \rm Consider the well-known Barenblatt solutions (recall that we work for $N \ge 3$) 
\[
u_B(\mathbf{x},t)=t^{-\lambda}\left(C-k\frac{|\mathbf{x}|^2}{t^{2\gamma}}\right)^{\frac1{m-1}}_+ \, , \ \lambda=\frac{N}{N(m-1)+2} \, ,  \ \gamma=\frac\lambda N \, , \ k=\frac{\lambda(m-1)}{2mN} \, ,
\]
which solve the Neumann problem \eqref{eq: pmeNeu} for sufficiently short times, provided $\mathbf{x}=\mathbf{0}$ belongs to $\Omega$ (else one just translates them). By a \emph{scaling invariant} estimate for solutions to \eqref{eq: pmeNeu} we mean a bound of the type
$$ \left\| u(t) \right\|_\infty  \le S(t,\left\| u_0 \right\|_1) \ \ \ \forall t>0 \, , \ \forall u_0 \in L^1(\Omega) $$
with $S(\cdot,\cdot)$ such that
$$ S(t,\left\| u_0 \right\|_1) = \left\| u_0 \right\|_1  S(\|u_0 \|_1^{m-1} \, t,1 )  \ \ \ \forall t>0 \, , \ \forall u_0 \in L^1(\Omega)  \, . $$
So $S(\cdot,\cdot)$ is completely determined by the one-variable function $S(y,1)$. Estimate \eqref{eq: soboStimaReg} is clearly scaling invariant. Moreover, it is a simple calculation to check that the left and right hand sides of \eqref{eq: soboStimaReg}, when evaluated on the Barenblatt solutions and considered e.g.\ in the time interval $[t,2t]$, have the same rate of divergence as $t\downarrow0$, and that the $L^1$ norm of $u_B(\cdot,t)$ is preserved along the evolution. This means that \eqref{eq: soboStimaReg} (for $q_0=1$) is sharp in the sense that there cannot hold another scaling invariant estimate with a better rate for $S(y,1)$ as $y \downarrow 0$. 
Also, we shall see in section \ref{sec: stimeWei} that, should the bound
\begin{equation}\label{eq: theta}\left\| u(t) \right\|_{\infty} \leq    K \left(t^{-\frac{\theta}{m-1}} \left\|  u_0 \right\|_{{q_0}^{\phantom{a}}}^{1-\theta} +  \left\| u_0 \right\|_{q_0}  \right)  \ \ \ \forall t>0 \, , \ \forall u_0\in L^{q_0}(\Omega)
\end{equation}
be valid for some $\theta<\frac{N(m-1)}{2q_0+N(m-1)}$ and for all $q_0\ge 1$ (in fact, $q_0\in[m, m+1)$ would suffice), then an embedding of $W^{1,2}(\Omega)$ into $L^q(\Omega)$ with $q>2N/(N-2)$ would be true as a consequence, a fact which is clearly false since $N\ge3$. 
\end{oss}

\begin{oss}\rm
We point out that we have stated our estimates for all $t>0$ rather than for \emph{almost every} $t>0$. This remains correct since in \cite[Ths. 11.2, 11.3]{Vaz07} continuity of solutions to \eqref{eq: pmeNeu} in $C([0,T];L^1(\Omega))$ (for all $T>0$) was proved, while in \cite{GMP} (see Section 3) the same was shown to hold for solutions to \eqref{eq: pmeNeuWei} at least when $\nu(\Omega)<\infty$ (which here will be the case), of course upon replacing $L^1(\Omega)$ with $L^1(\Omega;\nu)$. This comment also applies to the results of the forthcoming sections.
\end{oss}
\end{section}

\begin{section}{$L^\infty$ bounds for large time}\label{sec: stimeAsi}
In this section we shall prove some sharp improvements to the asymptotic estimates provided by Corollaries 1.3 and 1.4 of \cite{BG05m}. We remark that both the $L^{q_0}$-$L^\infty$ regularity estimates \eqref{eq: soboStimaReg} and \eqref{eq: stimaBG05} give no significant information about the behaviour of the solution $u(\cdot,t)$ as $t\rightarrow \infty$. So in order to obtain such informations one has to proceed through different techniques, whose starting point is the distinction between zero mean and non-zero mean solutions.

\subsection{The case $\mathbf{\overline{u}\boldsymbol{=}0}$}\label{sec: asiMn}
The following result, as we shall discuss in greater detail later, provides an improvement of the asymptotic estimates given by \cite[Cor. 1.3]{BG05m} for zero-mean solutions.
\begin{thm}\label{cor: soboMNAbsb}
Let $u_0 \in L^{q_0}(\Omega)$, with $q_0 \in [1,\infty)$ and $\overline{u_0}=0$. For the solution $u$ of \eqref{eq: pmeNeu} with initial datum $u_0$ the following estimate holds:
\begin{equation}\label{eq: soboAbsbMNullaPme}
\left\| u(t) \right\|_{\infty} \leq Q_1 \, t^{-\frac{N}{2q_0+N(m-1)}} \, \frac{1}{\left(Q_2 \, t + \left\| u_0 \right\|_{q_0}^{1-m}  \right)^{\frac{2 q_0}{(m-1)(2q_0+N(m-1))} } }     \ \ \ \forall t>0 \, ,
\end{equation}
where $Q_1$ and $Q_2$ are constants depending on $q_0$, $m$, $\Omega$. In particular, the absolute bound
\begin{equation}\label{eq: absB}
\left\| u(t)  \right\|_{\infty} \leq Q_3 \, t^{-\frac{1}{m-1}} \ \ \ \forall t>0
\end{equation}
holds true.
\begin{proof}
We start proving the assertion for $q_0>1$. In this case, the bound
\begin{equation}\label{eq: primoBound}
\left\| u(t) \right\|_{q_0} \leq \frac{1}{\left(B \, t + \left\| u_0 \right\|_{q_0}^{1-m}  \right)^{\frac{1}{m-1} } }
\end{equation}
($B=B(q_0,m,\Omega)$) was already established in the proof of Corollary 1.3 of \cite{BG05m}. Actually, this was done only for $q_0 \geq 2$, however one sees easily that the same result can be obtained also for $q_0 \in (1,2)$.

Now, let us consider estimate \eqref{eq: soboStimaReg} with the time origin shifted to $t/2$:
\begin{equation*}\label{eq: soboStimaRegShift}
\left\| u(t) \right\|_{\infty} \leq    K \left(\left( \frac{t}{2} \right)^{-\frac{N}{2q_0+N(m-1)}} \left\|  u(t/2) \right\|_{{q_0}^{\phantom{a}}}^{\frac{2q_0}{2q_0+N(m-1)}} +  \left\| u(t/2) \right\|_{q_0}  \right)  \, .
\end{equation*}
Exploiting \eqref{eq: primoBound} evaluated at time $t/2$ we remain with
\begin{equation*}\label{eq: soboStimaRegAsi}
\left\| u(t) \right\|_\infty \leq    K \left( \left(  \frac{t}{2} \right)^{-\frac{N}{2q_0+N(m-1)}} \frac{1}{\left( B  \frac{t}{2} + \left\| u_0 \right\|_{q_0}^{1-m}  \right)^{\frac{2 q_0}{(m-1)(2q_0+N(m-1))} }  } +  \frac{1}{\left(B  \frac{t}{2} + \left\| u_0 \right\|_{q_0}^{1-m}  \right)^{\frac{1}{m-1} } }  \right)  ,
\end{equation*}
from which \eqref{eq: soboAbsbMNullaPme} follows immediately. The absolute bound \eqref{eq: absB} is a straightforward consequence of \eqref{eq: soboAbsbMNullaPme}. Notice that performing the above calculations in full detail one gets that $Q_1\to \infty$ and $Q_2\to0$ as $q_0\downarrow1$, thus making it impossible to perform a limit to handle the case $q_0=1$. So in order to prove \eqref{eq: soboAbsbMNullaPme} for $q_0=1$ we need a different argument. First we notice that the validity of such inequality for $q_0=2$ implies that
\[
\left\|u(t)\right\|_\infty\le Q_1\,t^{-\frac{\gamma_0}{m-1}}\|u_0\|^{1-\gamma_0}_2 \, , \ \ \ \gamma_0=\frac{N(m-1)}{4+N(m-1)}\,;
\]
hence we have also
\begin{equation}\label{part-reg}
\left\|u(t)\right\|_\infty\le Q_1\left(\frac t2\right)^{-\frac{\gamma_0}{m-1}}\left\|u(t/2)\right\|_{\infty^{\phantom{a}}}^{\frac{1-\gamma_0}2}\left\|u_0\right\|_1^{\frac{1-\gamma_0}2} \, .
\end{equation}
Using the strategy outlined in the proof of Theorem \ref{teo: soboTeoReg}, see in particular the calculations following \eqref{eq: soboStimaRegCaso3}, one can show that the validity of \eqref{part-reg} implies in turn the validity of
\[
\left\|u(t)\right\|_\infty\le Q_3\,t^{-\frac{\omega_0}{m-1}}\left\|u_0\right\|^{1-\omega_0}_1,\ \ \ \omega_0=\frac{N(m-1)}{2+N(m-1)}\,.
\]
Therefore it follows that
\begin{equation}\label{eq: stimaGrezza}
\left\|u(t)\right\|_2\le \left\| u(t) \right\|_{\infty^{\phantom{a}}}^{\frac 12} \!\! \left\| u_0 \right\|_1^{\frac 12}\le Q_3^{\frac12} \, t^{-\frac1{m-1}\frac{\omega_0}2}\|u_0\|^{1-\frac{\omega_0}2}_1 \, .
\end{equation}
Exploiting \eqref{eq: soboAbsbMNullaPme} on the time interval $[t/2,t]$ with the choice $q_0=2$ and \eqref{eq: stimaGrezza} at time $t/2$ yields
\[
\begin{split}
\left\| u(t) \right\|_{\infty} &\leq Q_1 \, \left(\frac t2\right)^{-\frac{\gamma_0}{m-1}} \, \frac{1}{\left(Q_2 \, \frac t2 + \left\| u\left(t/2\right) \right\|_{2}^{1-m}  \right)^{\frac{1-\gamma_0}{m-1} } } \le \\
&\le Q_1 \, \left(\frac t2\right)^{-\frac{\gamma_0}{m-1}} \, \frac{1}{\left(Q_2 \, \frac t2 + Q_3^{\frac{1-m}2} \left(\frac t2\right)^{\frac{\omega_0}2} \|u_0\|^{(1-m)\left(1-\frac{\omega_0}2\right)}_1\right)^{\frac{1-\gamma_0}{m-1} } } = \\
&= Q_1 \, \left(\frac t2\right)^{-\frac{\omega_0}{m-1}}\, \frac{1}{\left(\frac t2\right)^{\frac{\gamma_0-\omega_0}{m-1}}\left(Q_2 \, \frac t2 + Q_3^{\frac{1-m}2} \left(\frac t2\right)^{\frac{\omega_0}2} \|u_0\|^{(1-m)\left(1-\frac{\omega_0}2\right)}_1\right)^{\frac{1-\gamma_0}{m-1} } } = \\
&=Q_1 \, \left(\frac t2\right)^{-\frac{\omega_0}{m-1}}\, \frac{1}{\left(Q_2 \, \left(\frac t2\right)^{1-\frac{\omega_0}2} + Q_3^{\frac{1-m}2} \|u_0\|^{(1-m)\left(1-\frac{\omega_0}2\right)}_1\right)^{\frac{1-\gamma_0}{m-1} } } \, .
\end{split}
\]
Hence \eqref{eq: soboAbsbMNullaPme} follows also for $q_0=1$ upon redefining the multiplicative constants and using the fact that for any given $\delta \in (0,1)$ one has $2^{\delta-1}\left( a^\delta+b^\delta\right)\le (a+b)^\delta\le  a^\delta+b^\delta $ for all $a,b>0$ ($\delta=1-{\omega_0}/{2}$ to our purposes). \normalcolor
\end{proof}
\end{thm}
\begin{oss}\label{oss: mn} \rm  In \cite[Cor. 1.3]{BG05m} the following estimate was proved:
\begin{equation}\label{eq: stimaAsBG}
\left\| u(t) \right\|_\infty \leq \frac{C}{\left(B(t-1) + \left\| u_0 \right\|_{q_0}^{1-m} \right)^{\frac{1-\alpha}{m-1}} } \ \ \ \forall t>1 \, ,
\end{equation}
for suitable constants $B,C>0$, where $\alpha\in(0,1)$ is defined as in \eqref{eq: stimaBG05alph}. Clearly \eqref{eq: stimaAsBG} is weaker than \eqref{eq: soboAbsbMNullaPme}, since it provides us with a \emph{slower} decay rate for $\| u(t) \|_\infty$ as $t \rightarrow \infty$. 
Actually, a rate of order $t^{-1/(m-1)}$ for $t$ large was already known to hold for bounded initial data since the pioneering work \cite{AR81} (see Theorems 3.1 and 4.1 there), even though the estimate given there still depended in a nontrivial way on the $L^\infty$ norm of the initial datum and in particular was not extendible to data not belonging to $L^\infty(\Omega)$. The sharpness of such rate follows by considering the separable variable solutions constructed in \cite[Sect. 2]{AR81}.

Finally, notice that \eqref{eq: soboAbsbMNullaPme} for small $t$ is equivalent to \eqref{eq: soboStimaReg}.
\end{oss}

In the statement of Theorem \ref{cor: soboMNAbsb} we pointed out that the constants $Q_1$, $Q_2$ depend on the whole $\Omega$ (in which we then absorb their dependence on $N$ and $C_S$ as well). The reason traces back to the constant appearing in Lemma 3.2 of \cite{AR81} (exploited in the proof of \cite[Cor. 1.3]{BG05m}), from which $B$ in \eqref{eq: primoBound} depends, whose existence is established by means of an argument by contradiction and so whose relation with significant features of the domain (such as its Poincar\'e constant) is completely unknown. We recall that Lemma 3.2 of \cite{AR81} states that there exists a constant $C_0>0$ depending on $r>1$ and $\Omega$ such that for all integrable functions $\xi$ such that $\overline{\xi}=0$ and $\xi^r \in W^{1,2}(\Omega)$ the following inequality holds:
\begin{equation*}\label{eq: lemmaAR}
\left\| \xi^r \right\|_2 \leq C_0 \left\| \nabla{(\xi^r)} \right\|_2 \, .
\end{equation*}
In fact, the result is also true for $r\geq {1}/{2}$. In \cite[Lem. 5.8]{GMP} a generalization of the above lemma is given, although still obtained through an argument by contradiction.

\subsection{The case $\mathbf{\overline{u}\boldsymbol{\neq}0}$}\label{sec: asiMNon}
We now show an improvement of the asymptotic estimate provided by \cite[Cor. 1.4]{BG05m} for solutions with non-zero mean. When making use of the Poincar\'e inequality \eqref{eq: poinIntroMedWei} we shall implicitly assume that the constant $C_P$ is the best one for which such inequality holds: in other words $1/{C_P^2}=\lambda_1$, where $\lambda_1$ is the first non-zero eigenvalue of (minus) the Laplacian operator with Neumann boundary conditions.
\begin{thm}\label{eq: teoAsiMNon}
Let $u$ be the solution of \eqref{eq: pmeNeu} corresponding to an initial datum $u_0 \in L^1(\Omega)$ with $\overline{u_0}=\overline{u}\neq 0$. There exists a constant $G=G(m,C_S,N,|\Omega|,u)>0$ such that the following estimate holds:
 \begin{equation}\label{eq: StimaExpFin}
\left\| u(t)-\overline{u} \right\|_\infty \leq G \, e^{-\frac{m}{C_P^2}|\overline{u}|^{m-1} \, t   } \ \ \ \forall t \geq 1  \, .
\end{equation}
\begin{proof}
We start considering data in $L^2(\Omega)$. From the uniform convergence to the mean value proved in \cite[Cor. 1.4]{BG05m}, we deduce in particular that there exists a time $\tau_0=\tau_0(u)$ such that
\begin{equation}\label{eq: convUnifNonW}
\inf_{\mathbf{x} \in \Omega}{ |u(\mathbf{x},t)|}  \geq \frac{|\overline{u}|}{2} \ \ \ \forall t \ge \tau_0 \, .
\end{equation}
Using \eqref{eq: convUnifNonW} in the differential inequality solved by $\| u(s) -\overline{u} \|_{\varrho}^\varrho$ (for a fixed $\varrho>1$), one gets
\begin{equation}\label{eq: convUnifNonW2}
\frac{\mathrm{d}}{\mathrm{d}s} \left\| u(s) -\overline{u} \right\|_{\varrho}^\varrho \leq  -4 \frac{\varrho-1}{\varrho} m {\left[\frac{|\overline{u}|}{2}\right]}^{m-1}  \left\| \nabla \left(  | u(s)- \overline{u} |^{\frac{\varrho}{2}}  \right) \right\|_{2}^2 \ \ \ \forall s \geq \tau_0 \, .
\end{equation}
In the case $\varrho=2$, \eqref{eq: convUnifNonW2} immediately leads to
\begin{equation}\label{eq: decL2}
\left\| u(t) -\overline{u} \right\|_{2} \leq e^{- \frac{m }{C_P^2} {\left[\frac{|\overline{u}|}{2}\right]}^{m-1} (t-\tau_0) } \left\| u(\tau_0) - \overline{u} \right\|_{2} \ \ \ \forall t \ge \tau_0  \, ,
\end{equation} 
upon applying \eqref{eq: poinIntroMedWei} to the right hand side. Now we observe that, setting $t_k=\tau_0+(t-\tau_0)(1-2^{-k})$ and $p_{k+1}=\frac{N}{N-2} p_k$ (let $p_0=2$), integrating inequality \eqref{eq: convUnifNonW2} between $t_k$ and $t_{k+1}$ with $\varrho=p_k$ and proceeding with similar (but easier) computations as in the proof of Theorem \ref{teo: soboTeoReg} one arrives at the following estimate:
\begin{equation}\label{eq: convUnifNonW3}
\left\| u(t) - \overline{u} \right\|_{\infty} \leq D \left( \left( {\left[\frac{|\overline{u}|}{2}\right]}^{m-1}  (t-\tau_0)\right)^{-1} + 1  \right)^{\frac{N}{4}} \left\|u(\tau_0) - \overline{u} \right\|_{2} \ \ \ \forall t > \tau_0  \, ,
\end{equation}
for a suitable constant $D=D(m,C_S,N,|\Omega|)$. By replacing $\tau_0$ with $\frac{t+\tau_0}{2}$ in \eqref{eq: convUnifNonW3}, exploiting \eqref{eq: decL2} at time $\frac{t+\tau_0}{2}$ (instead of $t$) and using the non-expansivity of $\| u(s) -\overline{u}\|_2$ (consequence of \eqref{eq: convUnifNonW2}) we obtain
\begin{equation}\label{eq: stimaConvUnif2}
\left\| u(t) - \overline{u} \right\|_{\infty} \leq  G \, e^{-M t } \left\| u_0 - \overline{u} \right\|_{2} \ \ \ \forall t \geq \tau_0 +1 \, ,
\end{equation}
for another suitable positive constants $G=G(m,C_S,N,|\Omega|,u)$ and $M=M(m,C_P,|\overline{u}|)$. Let us take $\tau_1\geq \tau_0+1$ such that
\begin{equation}\label{eq: tauPrimo}
 \frac{1}{2} |\overline{u}| \leq  |\overline{u}| - G \, e^{-M t } \left\| u_0 - \overline{u} \right\|_{2}  \ \ \ \forall t \geq \tau_1 \, .
\end{equation}
Since
$$ |u(\mathbf{x},s)| \geq |\overline{u}|-|u(\mathbf{x},s)-\overline{u}| \geq |\overline{u}|- G \, e^{-M s } \left\| u_0 - \overline{u} \right\|_{2} \ \ \ \forall s \geq \tau_1 \, , \ \hbox{for a.e. } \mathbf{x} \in \Omega \, , $$
we have:
\begin{equation}\label{eq: convUnifxx}
\begin{split}
 \frac{\mathrm{d}}{\mathrm{d}s} \left\| u(s) -\overline{u} \right\|_{2}^2 & =  -2m  \int_\Omega |u(\mathbf{x},s)|^{m-1}  \left| \nabla u(\mathbf{x},s) \right|^2  \, \mathrm{d}\mathbf{x} \leq \\
& \leq -2m \left( |\overline{u}|- G \, e^{-M s } \left\| u_0 - \overline{u} \right\|_{2} \right)^{m-1} \left\| \nabla{u}(s) \right\|_2^2 \leq \\
& \leq  - \frac{2m}{C_P^2} \left( |\overline{u}|- G \, e^{-M s } \left\| u_0 - \overline{u} \right\|_{2} \right)^{m-1}  \left\| u(s)-\overline{u} \right\|_2^2  \ \ \ \forall s \geq \tau_1 \, .
\end{split}
\end{equation}
By integrating \eqref{eq: convUnifxx} from $\tau_1 $ to $t$ (and again exploiting the non-expansivity of $\|u(s)-\overline{u} \|_2$) we get:
\begin{equation}\label{eq: preStimaExp}
\left\| u(t)-\overline{u} \right\|_2^2 \leq \left\| u_0-\overline{u} \right\|_2^2 \, e^{-\frac{2m}{C_P^2} \, \int_{\tau_1}^t \left( |\overline{u}|- G \, e^{-M s } \left\| u_0 - \overline{u} \right\|_{2} \right)^{m-1} \, \mathrm{d}s} \, .
\end{equation}
Setting $\epsilon= G \, e^{-M s } \left\| u_0 - \overline{u} \right\|_{2}$, thanks to the elementary inequalities (recall \eqref{eq: tauPrimo})
$$ \left(|\overline{u}|-\epsilon \right)^{m-1} \geq |\overline{u}|^{m-1}-(m-1)|\overline{u}|^{m-2} \epsilon \ \ \ \forall{m \geq 2} \, ,$$
$$ \left(|\overline{u}|-\epsilon \right)^{m-1} \geq |\overline{u}|^{m-1}-(m-1)(|\overline{u}|-\epsilon)^{m-2} \epsilon \geq |\overline{u}|^{m-1}-2^{2-m}(m-1)|\overline{u}|^{m-2} \epsilon  \ \ \ \forall{m \in (1,2) } \, ,$$
and to a straightforward computation of the time integral, estimate \eqref{eq: preStimaExp} reads
\begin{equation}\label{eq: preStimaExp2}
\left\| u(t)-\overline{u} \right\|_2 \leq G \, e^{-\frac{m}{C_P^2}|\overline{u}|^{m-1} \, t  } \ \ \ \forall t\geq \tau_1 \, ,
\end{equation}
being $G$ another suitable constant, different from the one above but depending on the same quantities (we can absorb the Poincar\'e constant $C_P$ into $C_S$). Combining \eqref{eq: stimaConvUnif2} with the time origin shifted to $t-\tau_0 -1$ and \eqref{eq: preStimaExp2} evaluated at time $t-\tau_0-1 $ (in place of $t$) and setting $\tau=\tau_0+\tau_1+1$ we finally arrive at
\begin{equation}\label{eq: preStimaExpFin}
\left\| u(t)-\overline{u} \right\|_\infty \leq G \, e^{-\frac{m}{C_P^2}|\overline{u}|^{m-1} \, t  } \ \ \ \forall t \geq \tau \, ,
\end{equation}
again for some $G=G(m,C_S,N,|\Omega|,u)$. The passage to general $L^1$ data and the fact that one can take, up to a different constant $G$, $\tau=1$ (or any other fixed $\tau>0$) in \eqref{eq: preStimaExpFin} follow from the $L^1$-$L^\infty$ regularizing effect.
\end{proof}
\end{thm}

\begin{oss}\label{oss: mdiff2}\rm
Estimate \eqref{eq: StimaExpFin} improves the one given in \cite[Cor. 1.4]{BG05m} with respect to time dependence. Indeed in \cite{BG05m} the authors provided a decay rate of order
$$G \, e^{-\frac{(1-\alpha)m}{C_P^2 2^{m-1}}|\overline{u}|^{m-1} t } $$
for $\|u(t)-u_0 \|_\infty $ ($\alpha$ is still defined by \eqref{eq: stimaBG05alph}), which is clearly slower than the one obtained here. However, we remark that an exponential rate of the type
$$G_\varepsilon \, e^{-\varepsilon \frac{m}{C_P^2}|\overline{u}|^{m-1} t } \ \ \ \forall \varepsilon \in (0,1) $$
could have also been obtained by exploiting the uniform convergence to the mean value proved in \cite{BG05m} together with \cite[Th. 3.3]{AR81}, although the constant $G_\varepsilon$ would diverge as $\varepsilon\uparrow1$. The new contribution of our Theorem \ref{eq: teoAsiMNon} lies in the fact that one can actually choose $\varepsilon=1$. Notice that writing $u(t)=\overline{u}+f(t)$ and linearizing \eqref{eq: pmeNeu} around the constant solution $\overline{u}$ one formally gets the equation $f_t=m|\overline{u}|^{m-1}\Delta f$, so that the rate given in Theorem \ref{eq: teoAsiMNon} coincides with the one predicted by linearization.
\end{oss}
In the following we shall prove that estimate \eqref{eq: StimaExpFin} is sharp by providing data for which a lower bound matching the upper bound in \eqref{eq: StimaExpFin} holds. The result we shall give depends on the well-known fact that, under our assumptions, the operator $-\Delta_{\mathcal{N}}$ (non-weighted) possesses a sequence of eigenvalues $\{\lambda_n\}$ (possibly repeated) and a corresponding sequence of eigenfunctions $\{\phi_n \}$ which form an orthonormal basis for $L^2(\Omega)$. 
\begin{pro}\label{lower}
Let $\Omega$ be a bounded $C^\infty$ domain.\ Consider the first non-zero eigenvalue $\lambda_1$ of (minus) the Laplacian with homogeneous Neumann boundary conditions, and the corresponding set of eigenfunctions (possibly consisting of one single element) normalized in $L^2(\Omega)$. Choose any of such eigenfunctions and denote it by $\psi_1$. Given $\overline{u_0}\neq 0$, let $u$ be the solution of \eqref{eq: pmeNeu} corresponding to the initial datum $u_0=\overline{u_0}+c_1 \psi_1$, with $|c_1|$ sufficiently small. Then there exists a constant $H=H(m,\Omega,c_1)>0$ such that the following estimate holds:
 \begin{equation}\label{eq: StimaExpFin2}
\left\| u(t)-\overline{u} \right\|_\infty \geq H \, e^{-\frac{m}{C_P^2}|\overline{u}|^{m-1} \, t   } \ \ \ \forall t > 0 \, .
\end{equation}
\begin{proof}
Let us consider first the solution to \eqref{eq: pmeNeu} corresponding to a regular initial datum $u_0$ such that $\overline{u_0}=1$ and $\| u_0 - 1 \|_\infty \leq \delta<1$, with $\delta$ to be chosen sufficiently small afterwards. Since $\| u(t)-1 \|_\infty$ is non-increasing (indeed $\| u(t)-1  \|_{\varrho}$ is non-increasing for all $\varrho\geq 1$, and so for $\varrho=\infty$ as well) we have that $|u(\mathbf{x},t)| \geq  1-\delta $ for all $\mathbf{x} \in \Omega$ and $t >0$. Hence by inserting such bound into the first line of \eqref{eq: convUnifxx} and solving similarly the resulting differential inequality we obtain the following estimate:
\begin{equation}\label{eq: bound1}
\left\| f(t)  \right\|_2 \leq \left\| f_0 \right\|_2 \, e^{-\lambda_1 m (1-\delta)^{m-1} t} \ \ \ \forall t>0 \, ,
\end{equation}
where we set $f(\mathbf{x},t)=u(\mathbf{x},t)-1$ and $f_0(\mathbf{x})=u_0(\mathbf{x})-1$. From standard quasilinear theory (see \cite{LSU}) we have that $f \in C^{\infty}(\overline{\Omega} \times (0,\infty))$. Also, straightforward computations show that $f$ solves the following differential equation:
\begin{equation}\label{eq: eqF}
f_t=m\Delta{f} + F \, ,
\end{equation}
up to defining $F$ as
\begin{equation*}
F=m(m-1)(1+f)^{m-2}|\nabla{f}|^2 + m\left[(1+f)^{m-1}-1\right]\Delta{f} \, .
\end{equation*}
If we choose as initial datum $f_0(\mathbf{x})=c_1 \, \psi_1(\mathbf{x})$, under the condition
\begin{equation}\label{eq: c1varia}
|c_1| \leq \frac{\delta}{\| \psi_1 \|_\infty} \, ,
\end{equation}
we see that \eqref{eq: bound1} becomes
\begin{equation}\label{eq: bound2}
\left\| f(t)  \right\|_2 \leq |c_1| \, e^{-\lambda_1 m (1-\delta)^{m-1} t} \ \ \ \forall t>0 \, .
\end{equation}
Define as usual $|f|_{C^0(\overline{\Omega})}:=\|f\|_\infty$ and, for any multi-index $\eta=(\eta_1,\ldots,\eta_N)$, the quantity $|\eta|=\eta_1+\ldots+\eta_N$ and the seminorms
\[
\left|f\right|_{C^k(\overline{\Omega})}:=\max_{|\eta|=k}\left\|\partial^\eta f\right\|_\infty \, ,\ \ \ k\in{\mathbb N} \, .
\]
From the uniform parabolicity of the equation at hand we infer the existence of a constant $Q=Q(\delta)>0$ such that
\begin{equation}\label{eq: stimaRegLady}
\left| f(t) \right|_{C^k(\overline{\Omega})} \leq Q \ \ \ \forall t>0 \, , \ \forall k \in \mathbb{N} \, ,
\end{equation}
where $Q$ can be taken to be independent of $c_1$ subject to \eqref{eq: c1varia}. We now recall the generalized interpolation inequalities
\begin{equation}\label{interpolation}
\left|g \right|^{\phantom{\frac 12}}_{C^j(\overline{\Omega})}\le C_{j,k,p}\, \left|g\right|_{C^k(\overline{\Omega})}^{\frac{N+jp}{N+kp}}\,\left\|g\right\|_{p^{\phantom{a}}}^{\frac{p(k-j)}{N+kp}} \, ,
\end{equation}
valid for all integers $k>j\ge0$ and real $p\ge1$ (see \cite[p. 130]{N} or, for a short review, \cite[App. 3]{BGV}). Combining \eqref{eq: bound2}, \eqref{eq: stimaRegLady} and \eqref{interpolation} we deduce that the following bounds hold:
\begin{equation*}\label{eq: bound3}
\left| f(t) \right|_{C^{l}(\overline{\Omega})} \leq Q^\prime |c_1|^{1-\epsilon} \, e^{-\lambda_1 m (1-\delta)^{m-1} (1-\epsilon) t} \ \ \ {\rm for}\ l=0,1,2  \, , \ \forall t>0 \, , \ \forall \epsilon \in (0,1)  \, ,
\end{equation*}
being $Q^\prime(\delta,\epsilon)>0$ another constant independent of $c_1$ subject to \eqref{eq: c1varia}. From the very definition of $F$ and from the fact that the $L^\infty$ norm of $f(t,\cdot)$ is non-increasing we have that, for suitable constants $Q_0=Q_0(\delta)$, $Q^{\prime\prime}=Q^{\prime\prime}(\delta,\epsilon)$,
\begin{equation}\label{eq: boundG}\begin{split}
\left\| F(t) \right\|_{\infty} &\leq Q_0\left(\left| f(t) \right|_{C^1(\overline{\Omega})}^2+ \left\| f(t) \right\|_{\infty} \, \left| f(t) \right|_{C^2(\overline{\Omega})}\right)\normalcolor \le \\
&\leq Q^{\prime\prime} |c_1|^{2(1-\epsilon)} \, e^{- 2 \lambda_1 m (1-\delta)^{m-1} (1-\epsilon) t} \ \ \ \forall t>0 \, , \forall \epsilon \in (0,1)  \, .
\end{split}\end{equation}
We want to study the asymptotic behaviour of $\alpha_1(t)=\langle f(t) , \psi_1 \rangle_{L^2(\Omega)} $, that is the first Fourier coefficient of $f$ along the evolution. 
In order to do that, let us multiply equation \eqref{eq: eqF} by $\psi_1$ and integrate on $\Omega$, so as to obtain the following differential equation for $\alpha_1(t)$:
\begin{equation*}\label{eq: eqC}
\dot{\alpha}_1(t)= - \lambda_1 m \alpha_1(t) + \left\langle F(t) , \psi_1 \right\rangle_{L^2(\Omega)} \, .
\end{equation*}
Duhamel principle entails that $\alpha_1$ must satisfy the following nonlinear integral equation (recall that by construction $\alpha_1(0)=c_1$): \normalcolor
\begin{equation}\label{eq: eqCSolved}
\alpha_1(t)= e^{-\lambda_1 m t} \left[ c_1 + \int_0^t e^{\lambda_1 m s} \left\langle F(s) , \psi_1 \right\rangle_{L^2(\Omega)} \, \mathrm{d}s \right] \, .
\end{equation}
If we choose $\delta$ and $\epsilon$ sufficiently small, namely such that $2(1-\delta)^{m-1} (1-\epsilon) > 1$, thanks to \eqref{eq: boundG} we easily see that the time integral in \eqref{eq: eqCSolved} can be bounded in the following way:
\begin{equation*}\label{eq: boundInteg}
 \left| \int_0^t e^{\lambda_1 m s} \left\langle F(s) , \psi_1 \right\rangle_{L^2(\Omega)} \, \mathrm{d}s \right| \leq B \, |c_1|^{2(1-\epsilon)} \, ,
\end{equation*}
where $B$ is a suitable constant independent of $c_1$ and $t\ge0$. Pick now $|c_1|$ small enough so that $|{c}_1|>B\,|{c}_1|^{2(1-\epsilon)}$, this being possible since $2(1-\epsilon)> 1$. Under such bound on $|c_1|$ we then deduce from \eqref{eq: eqCSolved} that the \emph{exact} decay rate of $\alpha_1(t)$ is $e^{-\lambda_1 m t}$, in the sense that $|\alpha_1(t)|\sim  e^{-\lambda_1 m t}$ as $t\to\infty$. Hence \eqref{eq: StimaExpFin2} holds true since there exists a suitable constant $\widetilde K>0$ such that
\[
\left\|u(t)-1\right\|_\infty=
\left\|f(t)\right\|_\infty\ge |\Omega|^{-\frac12}\left\|f(t)\right\|_2\ge|\Omega|^{-\frac12}\left|\alpha_1(t)\right|\ge \widetilde K \, e^{-\lambda_1 m t} \ \ \ \forall t>0 \, .
\]
The case of initial data with non-zero mean value $\overline{u_0}\neq 1$ can be brought back to the case $\overline{u_0}=1$ by means of a standard time scaling argument.
\end{proof}
\end{pro}
\end{section}

\begin{section}{The weighted case}\label{sec: stimeWei}
As already mentioned, the main results provided in the previous sections can be extended to the case of the weighted porous media equation with Neumann boundary conditions \eqref{eq: pmeNeuWei}, dealt with in detail in \cite{GMP}, under the sole hypotheses that the weights $\rho_\nu,\rho_\mu$ satisfy Sobolev-type inequalities like \eqref{eq: soboIntroMedWei} or the weaker \eqref{eq: sobDebWei}. Since the proofs of the corresponding theorems basically would be the same as for the non-weighted case, we shall only give statements with short comments, except in the cases of Theorem \ref{thm: impInv}, Theorem \ref{teo: convUnif} and related lemmas, which require new arguments. Before doing that, we briefly recall some notation and essential definitions concerning weighted Sobolev spaces.

Let $\Omega \subset \mathbb{R}^N$ be a domain and let $\nu$ and $\mu$ be two measures defined on it, both absolutely continuous with respect to the Lebesgue measure. We indicate as $\rho_{\nu}$ and $\rho_{\mu}$ the corresponding weights (or densities), which will always be assumed to be strictly positive almost everywhere. For all $p \in [1,\infty)$ we introduce the Banach space $L^p(\Omega;\nu)$ of equivalence classes of Lebesgue measurable functions $f$ such that
$$ \left\| f \right\|_{p;\nu}^p = \int_{\Omega} |f|^p \, \mathrm{d}\nu = \int_{\Omega} |f|^p \, \rho_{\nu}\mathrm{d}\mathbf{x} < \infty \, . $$
The same applies for the measure $\mu$. According to \cite{KO84}, we define the weighted Sobolev space $W^{1,p}(\Omega; \nu, \mu)$ as the space of all (equivalence classes of) functions $v \in W^{1,1}_{loc}(\Omega)$ such that
\begin{equation*}
\left\| v \right\|^p_{p;\nu,\mu} = \left\| v \right\|^p_{p;\nu} + \left\| \nabla v \right\|^p_{p;\mu} < \infty \, .
\end{equation*}
Without further assumptions on $\rho_{\nu}$ and $\rho_{\mu}$ in general $W^{1,p}(\Omega; \nu, \mu)$ would not be complete.
\begin{den} \label{def: Bp}
For all $p \in (1,\infty)$ we denote as $B^p(\Omega)$ the class of all measurable functions $g$ such that $g > 0$ a.e.\ and
\begin{equation*}
f^{\frac{1}{1-p}} \in L^1_{loc}(\Omega) \, .
\end{equation*}
\end{den}
One can prove \cite[Th. 2.1]{KO84} that if $p \in (1,\infty)$ and $\rho_{\mu} \in B^p(\Omega)$ then $W^{1,p}(\Omega; \nu, \mu)$ is indeed complete. If $p=1$ the same result holds true providing that the condition $\rho_{\mu} \in B^p(\Omega)$ is replaced by $\rho_{\mu}^{-1} \in L^{\infty}_{loc}(\Omega)$.

The fact that for any $\varphi \in C^{\infty}_c(\Omega)$ the quantity \mbox{$\| \varphi \|_{p;\nu,\mu} $} is finite is equivalent (see \cite[Lem. 4.4]{KO84}) to the local finiteness of $\nu$ and $\mu$, that is
\begin{equation}\label{eq: condizioneFin}
\rho_{\nu}, \rho_{\mu}  \in L^1_{loc}(\Omega)  \, .
\end{equation}
Assuming \eqref{eq: condizioneFin}, one can define the space $W^{1,p}_0(\Omega; \nu, \mu)$ as the closure of $C^{\infty}_c(\Omega)$ in $W^{1,p}(\Omega; \nu, \mu)$.

Now we are ready to state the weighted counterpart of the $L^{q_0}$-$L^\infty$ regularity Theorem \ref{teo: soboTeoReg}. Hereafter, the assumptions on the weights (inner regularity and boundedness away from zero) and the concepts of solution will always be the ones given in Section \ref{sec:def}.
\begin{thm}\label{teo: soboTeoRegWei}
Let $\nu(\Omega)<\infty$ and let inequality \eqref{eq: sobDebWei} hold true for some $\sigma>1$. Then for the solution $u$ of \eqref{eq: pmeNeuWei} corresponding to an initial datum $u_0 \in L^{q_0}(\Omega;\nu)$ with $q_0 \in [1,\infty)$ the following estimate holds:
\begin{equation}\label{eq: soboStimaRegWei2}
\left\| u(t) \right\|_{\infty} \leq    K \left(t^{-\frac{\sigma}{(\sigma-1)q_0+\sigma(m-1)}} \left\|  u_0 \right\|_{q_0;\nu^{\phantom{a}}}^{\frac{(\sigma-1)q_0}{(\sigma-1)q_0+\sigma(m-1)}} +  \left\| u_0 \right\|_{q_0;\nu}  \right)  \ \ \ \forall t>0 \, ,
\end{equation}
where $K$ is a constant which depends only on $m$, $C_S$, $\sigma$ and $\nu(\Omega)$.
\end{thm}
The proof of Theorem \ref{teo: soboTeoRegWei} is exactly the same as the one of Theorem \ref{teo: soboTeoReg}: indeed, going back to the latter, one easily sees that in order to repeat it the only relevant assumption is the fact that the weaker inequality \eqref{eq: sobDebWei} holds (with $\sigma=N/(N-2)$ in that case).

Now we show that the $L^{q_0}$-$L^\infty$ regularity estimate \eqref{eq: soboStimaRegWei2} is in fact the correct one, namely its validity implies in turn the validity of a Sobolev embedding for $W^{1,2}(\Omega;\nu,\mu)$.
\begin{thm}\label{thm: impInv}
Let $\nu(\Omega)<\infty$ and suppose that there exist a constant $K>0$ and two given numbers $\sigma>1$ and $q_0 \in [m,m+1)$ such that, for all $u_0\in L^{q_0}(\Omega;\nu)$, the solution $u$ of \eqref{eq: pmeNeuWei} corresponding to the initial datum $u_0$ satisfies the following estimate:
\begin{equation}\label{eq: soboStimaRegWeiInv}
\left\| u(t) \right\|_{\infty} \leq    K \left(t^{-\frac{\sigma}{(\sigma-1)q_0+\sigma(m-1)}} \left\|  u_0 \right\|_{q_0;\nu^{\phantom{a}}}^{\frac{(\sigma-1)q_0}{(\sigma-1)q_0+\sigma(m-1)}} +  \left\| u_0 \right\|_{q_0;\nu}  \right)  \ \ \ \forall t>0 \, .
\end{equation}
Then there exists a constant $C^\prime>0$ such that the functional inequality
\begin{equation}\label{eq: impInvNeuDF}
\left\|  v \right\|_{2\sigma;\nu} \leq C^\prime \left( \left\| \nabla{v}  \right\|_{2;\mu} + \left\| v \right\|_{\frac{q_0}{m};\nu} \right) \ \ \ \forall v \in W^{1,2}(\Omega;\nu,\mu)
\end{equation}
holds as well.
\begin{proof}
We proceed along the lines of the proofs of Theorems 4.3 and 5.6 of \cite{GMP}. First of all, consider a non-negative initial datum $u_0 \in L^{\infty}(\Omega) \cap W^{1,2}(\Omega;\nu,\mu)$ and its corresponding solution to \eqref{eq: pmeNeuWei} $u$. Exploiting the classical interpolation inequality between the exponents $q_0 \in [m,m+1)$, $m+1$, $\infty$ and the non-expansivity of the $L^{q_0}(\Omega;\nu)$ norm, one easily obtains that
\begin{equation}\label{eq: interp}
\left\| u(t) \right\|_{{m+1};\nu} \leq \left\| u(t) \right\|_{\infty^{\phantom{a}}}^{\frac{m+1-q_0}{m+1}} \left\| u_0 \right\|_{q_0;\nu^{\phantom{a}}}^{\frac{q_0}{m+1}} \ \ \ \forall t>0 \, ;
\end{equation}
combining \eqref{eq: interp} with estimate \eqref{eq: soboStimaRegWeiInv} one gets, up to a multiplicative constant that for the sake of simplicity from now on we shall always denote as $K$,
\begin{equation}\label{eq: interpReg}
\left\| u(t) \right\|_{m+1;\nu}^{m+1} \leq K \left[ t^{-\frac{\sigma(m+1-q_0)}{(\sigma-1)q_0+\sigma(m-1)}} \left\| u_0 \right\|_{q_0;\nu^{\phantom{a}}}^{\frac{q_0(2\sigma m -m-1)}{(\sigma-1)q_0 +\sigma(m-1) }} + \left\| u_0 \right\|_{q_0;\nu^{\phantom{a}}}^{m+1}  \right] \ \ \ \forall t>0 \, .
\end{equation}
Now we need to use the following fundamental inequality:
\begin{equation}\label{eq: impInvPmeDir}
\left\| u(t) \right\|_{m+1;\nu}^{m+1} - \left\| u_0 \right\|_{m+1;\nu}^{m+1} \geq - (m+1) \, t \left\| \nabla{(u_0^m)}  \right\|_{2;\mu}^2  \ \ \ \forall t>0 \, ,
\end{equation}
which was proved in \cite{GMP} (see again Theorems 4.3 and 5.6). From \eqref{eq: interpReg} and \eqref{eq: impInvPmeDir} we deduce that
\begin{equation}\label{eq: preSobo1}
\begin{split}
& \left\| u_0 \right\|_{m+1;\nu}^{m+1} \leq \\
& \leq K \left[  t^{-\frac{\sigma(m+1-q_0)}{(\sigma-1)q_0+\sigma(m-1)}} \left\| u_0 \right\|_{q_0;\nu^{\phantom{a}}}^{\frac{q_0(2\sigma m -m-1)}{(\sigma-1)q_0 +\sigma(m-1) }} + (m+1) t \left\| \nabla{(u_0^m)}  \right\|_{2;\mu}^2  +\left\| u_0 \right\|_{q_0;\nu^{\phantom{a}}}^{m+1}  \right] \ \ \ \forall t>0 \, .
\end{split}
\end{equation}
Minimizing the right hand side of \eqref{eq: preSobo1} with respect to $t>0$ one arrives at
\begin{equation}\label{eq: preSobo2}
\left\| u_0 \right\|_{m+1;\nu}^{m+1} \leq K \left[ \left\| u_0 \right\|_{q_0;\nu^{\phantom{a}}}^{\frac{q_0(2\sigma m -m-1)}{2\sigma m -q_0 }} \left\| \nabla{(u_0^m)}  \right\|_{2;\mu}^{\frac{2\sigma(m+1-q_0)}{2\sigma m-q_0 }}  +\left\| u_0 \right\|_{q_0;\nu^{\phantom{a}}}^{m+1}  \right]  \, .
\end{equation}
Proceeding exactly as explained in the proof of Theorem 5.6 of \cite{GMP}, that is by approximating $u_0^{1/m}$ with a sequence of regular functions, one can show that \eqref{eq: preSobo2} is equivalent to the following inequality:
\begin{equation}\label{eq: preSobo3}
\left\| v \right\|_{r;\nu} \leq K \left[ \left(\left\| \nabla{v} \right\|_{2;\mu} + \left\| v \right\|_{\frac{q_0}{m};\nu} \right)^{\vartheta} \left\| v \right\|_{s;\nu}^{1-\vartheta} \right]  \ \ \forall v\geq 0  \in L^{\infty}(\Omega) \cap W^{1,2}(\Omega;\nu,\mu) , \ \frac{1}{r}=\frac{\vartheta}{q}+\frac{1-\vartheta}{s}  ,
\end{equation}
where it is understood that
$$ r=\frac{m+1}{m} \, , \ q=2\sigma \, , \ s=\frac{q_0}{m} \, , \ \vartheta=\frac{2\sigma m(m+1-q_0)}{(m+1)(2\sigma m -q_0)} \, .$$
The fundamental result \cite[Th. 3.1]{BCLS95} ensures that (up to a different multiplicative constant) one can put $\vartheta=1$ and thus $r=2\sigma$ in \eqref{eq: preSobo3}, this leading to the claimed inequality \eqref{eq: impInvNeuDF} at least for regular positive functions; the fact that such inequality holds in the whole space $W^{1,2}(\Omega;\nu,\mu)$ then follows by means of a density argument and by writing $v=v_{+}-v_{-}$.
\end{proof}
\end{thm}
A straightforward consequence of the above Theorems \ref{teo: soboTeoRegWei} and \ref{thm: impInv} is the following, which is one of the main results of the paper.
\begin{cor}\label{cor: equivNeu}
Let $\nu(\Omega)<\infty$. The validity of inequality \eqref{eq: sobDebWei} for some $\sigma>1$ is \emph{equivalent} to the validity of the family of estimates \eqref{eq: soboStimaRegWei2} for all $q_0 \in [1,\infty)$.
\end{cor}
\begin{oss}\rm
We shall provide a specific example in which the bound \eqref{eq: soboStimaRegWei2} (for $q_0=1$) is sharp in the sense that no better scaling invariant estimate can hold, as in the non-weighted case discussed in Remark \ref{oss: barenblatt}.

The setting is one-dimensional, in particular we choose $\Omega=(0,1)$. We consider the weights $\rho_\nu(x)=1$, $\rho_\mu(x)=x^\beta$. The Sobolev inequality \eqref{eq: soboIntroMedWei} associated to such weights is known to hold with $\sigma=\frac1{\beta-1}$ if $\beta\in(1,2)$ (see Section \ref{sec: exaSobo} for appropriate references). An explicit calculation shows that the Barenblatt-type functions
\begin{equation}\label{eq: barenblattbeta}
u_{B,\beta}(x,t)=t^{-\zeta}\left(C-k\frac{x^{\omega}}{t^{\omega\zeta}}\right)^{\frac1{m-1}}_+,
\end{equation}
where $C>0$ is a free (mass) paramenter, $\zeta=\frac1{m+1-\beta}$, $\omega=2-\beta$ and $k=\frac{m-1}{m(2-\beta)(m+1-\beta)}$, are solutions to the corresponding weighted equation for any $t>0$. We consider them for $t$ small enough to ensure that the support of $u_{B,\beta}(\cdot, t)$ is bounded away from the point $x=1$. It is straightforward to prove that the mass of $u_{B,\beta}(\cdot, t)$ is preserved in time, so that also the Neumann boundary condition at $x=0$ is satisfied. One then checks directly that such solutions belong to the appropriate functional spaces, so that they are energy solutions (for $t>0$) in the sense of Definition \ref{den: solDebNeu}.

We notice that (minus) the power of time appearing in the bound \eqref{eq: soboStimaRegWeiInv} coincides, for $q_0=1$ and $\sigma=\frac1{\beta-1}$, with $\frac1{m+1-\beta}$, namely with the power of time $\zeta$ appearing in \eqref{eq: barenblattbeta}. The fact that there cannot hold another scaling invariant estimate with a better rate then follows as in Remark \ref{oss: barenblatt}.
\end{oss}

We want to stress that Corollary \ref{cor: equivNeu} has a natural counterpart for the Dirichlet problem
\begin{equation} \label{eq: pmeDirWei}
\begin{cases}
u_t =\rho_{\nu}^{-1}  \operatorname{div}\left(\rho_{\mu} \, \nabla{\left( u^m \right)} \right) & \textnormal{in} \  \Omega\times(0,\infty)  \\
 u =0 & \textnormal{on} \ \partial\Omega \times (0,\infty) \\
 u(\cdot,0)=u_0(\cdot) & \textnormal{in} \ \Omega
 \end{cases} \, ,
\end{equation}
whose well-posedness was also analysed in \cite{GMP} and whose correct $L^{q_0}$-$L^\infty$ regularity estimates basically had already been proved in \cite[Th. 1.5]{BG05m} (there is no difficulty in extending the proof provided therein to our context). Here is the statement:
\begin{cor}\label{teo: poinTeoReg}
Let $\nu(\Omega)<\infty$. The validity of the Sobolev inequality
\begin{equation}\label{eq: soboDirWei}
\left\| v \right\|_{2\sigma;\nu} \leq C_S \left\| \nabla{v} \right\|_{2;\mu} \ \ \ \forall v \in W^{1,2}_0(\Omega;\nu,\mu)
\end{equation}
for some $\sigma>1$ and some constant $C_S>0$ is \emph{equivalent} to the validity of the family of estimates
\begin{equation}\label{eq: soboStimaRegDirWei}
\left\| u(t) \right\|_{\infty} \leq    K \, t^{-\frac{\sigma}{(\sigma-1)q_0+\sigma(m-1)}} \left\|  u_0 \right\|_{q_0;\nu^{\phantom{a}}}^{\frac{(\sigma-1)q_0}{(\sigma-1)q_0+\sigma(m-1)}}   \ \ \ \forall t>0 \, , \ \forall q_0\in[1,\infty) \, ,
\end{equation}
where $u$ is the solution of \eqref{eq: pmeDirWei} corresponding to a generic initial datum $u_0 \in L^{q_0}(\Omega;\nu)$ and $K>0$ is a suitable constant independent of $u_0$.
\end{cor}
The finiteness of the measure, in this case, is necessary only in order to prove the converse (i.e.\ the fact that \eqref{eq: soboStimaRegDirWei} implies \eqref{eq: soboDirWei}, see \cite[Th. 4.3]{GMP}).

As for the asymptotic behaviour of solutions, we begin with an analogue of Theorem \ref{cor: soboMNAbsb}.
\begin{thm}\label{cor: soboMNAbsbWei}
Let $\nu(\Omega)<\infty$ and let the Sobolev-type inequality \eqref{eq: soboIntroMedWei} hold true for some $\sigma>1$. Then if $u_0 \in L^{q_0}(\Omega;\nu)$, with $q_0 \in [1,\infty)$, $\overline{u_0}=0$ and $u$ is the solution of \eqref{eq: pmeNeuWei} corresponding to the initial datum $u_0$, the following estimate holds:
\begin{equation*}\label{eq: soboAbsbMNullaPmeWei}
\left\| u(t) \right\|_{\infty} \leq Q_1 \, t^{-\frac{\sigma}{(\sigma-1)q_0+\sigma(m-1)}} \, \frac{1}{\left(Q_2 t + \left\| u_0 \right\|_{q_0;\nu}^{1-m}  \right)^{\frac{(\sigma-1) q_0}{(m-1)[(\sigma-1)q_0+\sigma(m-1)]} } }  \ \ \ \forall t>0 \, ,
\end{equation*}
where $Q_1$ and $Q_2$ are constants depending on $q_0$, $m$, $\Omega$, $\nu$, $\mu$.
\end{thm}
Again, the proof of this result is similar to the one of Theorem \ref{cor: soboMNAbsb}. Actually, in order to be able to exploit Lemma 3.2 of \cite{AR81}, a priori one should require that the embedding $W^{1,2}(\Omega;\nu,\mu) \hookrightarrow L^2(\Omega;\nu) $ is compact. In principle this is not a major restriction since it is quite common that, at least on bounded domains, the validity of a Sobolev-type inequality comes together with the compactness of such embedding (see the examples of Section 6 of \cite{GMP}). Nonetheless, as previously mentioned, in \cite[Lem. 5.8]{GMP} a generalized version of that lemma is provided, which does not require compactness.

As concerns solutions to \eqref{eq: pmeNeuWei} with non-zero mean, in order to obtain the weighted version of the asymptotic estimate given in Theorem \ref{eq: teoAsiMNon} first we have to prove uniform convergence of such solutions to their mean value. Since the proof of \cite[Cor. 1.4]{BG05m} (which takes advantage of the uniform spatial H\"{o}lder continuity of solutions to \eqref{eq: pmeNeu}) in the generality of present context does not work, we shall provide a different (functional) argument, which basically only uses the assumed validity of the Sobolev-type inequality \eqref{eq: soboIntroMedWei}. To this end, we need two preliminary lemmas.

\begin{lem}\label{lem: ratio}
Given $r \geq 1/2 $, $m > 1$ and a fixed constant $R>1$, let
\begin{equation*}\label{eq: defPhi}
\Phi_{r,m}(x)=\int_0^x{|y|^{r-1}|y+1|^{\frac{m-1}{2}} \, \mathrm{d}y } \ \ \ \forall x \in [-R,R]  \, .
\end{equation*}
Then there exist two positive constants $C_1=C_1(m)$ and $C_2=C_2(m,R)$ such that
\begin{equation}\label{eq: lowUppBound}
\frac{C_1}{r^{1+ 1 \vee \left[\frac{m-1}{2}\right] }} |x|^r \leq |\Phi_{r,m}(x)| \leq \frac{C_2}{r} |x|^r \ \ \ \forall x \in [-R,R]  \, .
\end{equation}
\begin{proof}
One has to study the ratio $ |\Phi_{r,m}(x)|/|x|^r $. The bound from above is easily obtainable, since
$$ |\Phi_{r,m}(x)| \leq (R+1)^{\frac{m-1}{2}} \left| \int_0^x  |y|^{r-1}   \, \mathrm{d}y \right|  = \frac{(R+1)^{\frac{m-1}{2}}}{r} |x|^r \ \ \ \forall x \in [-R,R]  \, .  $$
In order to get a lower bound, we begin with the easier case $m=3$. Recall that, when $x\not=0$, we use the convention $x^r:=|x|^{r-1}x$. First of all note that
\begin{equation}\label{eq: ricercaMin1}
\Phi_{r,m}(x)/x^r \geq \frac{\int_0^x x^{r-1} \, \mathrm{d}y }{x^r} = \frac{1}{r}  \ \ \ \forall x \in (0,R] \, ;
\end{equation}
also,
\begin{equation}\label{eq: ricercaMin2}
\Phi_{r,3}(x)/x^r=\frac{ \int_0^{|x|} y^{r-1}(1-y) \, \mathrm{d}y }{|x|^r} = \frac{1}{r}-\frac{|x|}{r+1} \ge \frac{1}{r(r+1)}  \ \ \  \forall x \in [-1,0) \, .
\end{equation}
So we are left to study the minimum of $ \Phi_{r,3}(x)/x^r $ as $x$ varies in $[-R,-1)$. We have:
\begin{equation*}\label{eq: rm0}
\begin{split}
 f_r(x)=\Phi_{r,3}(x)/x^r & =\frac{ \int_0^{|x|} y^{r-1}|y-1| \, \mathrm{d}y }{|x|^r} = \frac{ \int_0^{1} y^{r-1}(1-y) \, \mathrm{d}y + \int_{1}^{|x|} y^{r-1}(y-1) \, \mathrm{d}y}{|x|^r} = \\
& = \frac{\frac{2}{r(r+1)} + \frac{|x|^{r+1}}{r+1}-\frac{|x|^r}{r} } {|x|^r} = \frac{2}{r(r+1)}|x|^{-r} + \frac{|x|}{r+1} - \frac{1}{r}   \ \ \  \forall x \in [-R,-1) \, .
\end{split}
\end{equation*}
Notice that
$$\frac{{\rm d}}{{\rm d}s}\left(\frac {2}{r(r+1)} s^{-r} +\frac s{r+1}-\frac1r\right) =-\frac{2}{r+1} s^{-r-1} + \frac{1}{r+1} \, ,  $$
whose zero is attained at $s_0=2^{\frac{1}{r+1}}$, so that
\begin{equation}\label{eq: ricercaMinMed}
f_r(x)\ge\frac{1}{r}\left( 2^{\frac{1}{r+1}} - 1 \right)\ \ \  \forall x \in [-R,-1)\, . 
\end{equation}
We finally observe that
\begin{equation}\label{eq: ricercaMin3}
\frac{1}{r}\left( 2^{\frac{1}{r+1}} - 1 \right) \ge \frac{\log{2}}{r(r+1)} \, .
\end{equation}
By collecting together \eqref{eq: ricercaMin1}, \eqref{eq: ricercaMin2}, \eqref{eq: ricercaMinMed} and \eqref{eq: ricercaMin3} one gets the lower bound in \eqref{eq: lowUppBound} for $m=3$.

Now, let us consider the case $m > 3$. Since the function $g(s)=|s|^{(m-1)/2} $ is convex, from Jensen's inequality we have:
\begin{equation*}
\frac{\left|\Phi_{r,m}(x)\right|}{|x|^r} = \frac1r \, \frac{\int_0^{|x|}{y^{r-1}|y-1|^{\frac{m-1}{2}} \, \mathrm{d}y }}{\int_0^{|x|} \, y^{r-1} \mathrm{d}y }  \geq  r^{\frac{m-3}{2} } \left( \frac{ \int_0^{|x|}{y^{r-1}|y-1| \, \mathrm{d}y   }} {|x|^r} \right)^{\frac{m-1}{2}}   \ \ \ \forall x \in [-R,0) \, ,
\end{equation*}
and we can bound the right hand side from below just by applying to it the above estimates for $m=3$, so as to obtain
\begin{equation*}
\frac{\left|\Phi_{r,m}(x)\right|}{|x|^r} \geq  \frac{C_1}{r^{1+\frac{m-1}{2}}} \ \ \ \forall x \in [-R,0)
\end{equation*}
for some constant $C_1=C_1(m)$, which together with \eqref{eq: ricercaMin1} gives the claimed lower bound in \eqref{eq: lowUppBound}. We remain to deal with the case $m \in (1,3)$, that is when the function $g(s)$ is no more convex. Again, since \eqref{eq: ricercaMin1} holds for all $m>1$, we can restrict ourselves to analyse the ratio $ |\Phi_{r,m}(x)|/|x|^r $ for $x \in [-R,0)$. Since $(m-1)/{2}=\alpha \in (0,1)$, straightforward calculations show that
\begin{equation*}\label{eq: maggElem}
|y-1|^\alpha \geq  |y^\alpha - 1| \ \ \ \forall y \in [0,R] \, ;
\end{equation*}
therefore,
\begin{equation}\label{eq: casoSotto}
\frac{\left|\Phi_{r,m}(x)\right|}{|x|^r} =  \frac{\int_0^{|x|}{y^{r-1}|y-1|^{\frac{m-1}{2}} \, \mathrm{d}y }}{|x|^r} \geq   \frac{ \left.\int_0^{|x|}{y^{r-1} |y^{\frac{m-1}{2}}-1 | \, \mathrm{d}y }\right. }{|x|^r}  \ \ \ \forall x \in [-R,0) \, .
\end{equation}
From now on, thanks to \eqref{eq: casoSotto}, by means of computations analogous to the ones developed in the case $m=3$ one can prove that also the right hand side of \eqref{eq: casoSotto} is bounded from below by a constant times $1/r^2$.
\end{proof}
\end{lem}
The next lemma is a crucial one. Basically it is a refinement of the already mentioned Lemma 5.8 of \cite{GMP} in the particular case of the functions $\Phi_{r,m}(\cdot)$, and so its proof is similar to the one of the latter.
\begin{lem}\label{lem: funct-inequality}
Suppose that there exists a constant $C_{PS}$ such that the Poincar\'e-Sobolev inequality
\begin{equation}\label{eq: poinSob}
\left\| v - \overline{v} \right\|_{2\sigma;\nu} \leq C_{PS} \left\| \nabla{v} \right\|_{2;\mu} \ \ \ \forall v \in W^{1,2}(\Omega;\nu,\mu)
\end{equation}
holds for some $\sigma \geq 1$. Then for any $m>1$ and $R>1$ there exists a constant $C_\ast=C_\ast(m,\Omega, \nu,\mu,R)$ such that
\begin{equation}\label{eq: poinSobPhi}
\left\| \Phi_{r,m}(\xi) \right\|_{2\sigma;\nu} \leq C_\ast \left\| \nabla{\Phi_{r,m}(\xi)} \right\|_{2;\mu}
\end{equation}
$$
\forall r \geq 1/2 \, , \ \forall \xi : \left\| \xi \right\|_\infty \leq R \, , \ \overline{\xi}=0 \, , \ \Phi_{r,m}(\xi) \in W^{1,2}(\Omega;\nu,\mu) \, .
$$
\begin{proof}
As in \cite[Lem. 5.8]{GMP} we proceed by contradiction. That is, if such a constant $C_\ast$ did not exist then there would be a sequence of numbers $r_n \geq 1/2$ and a corresponding sequence of non-identically-zero functions $\{\xi_n\}$ such that $r_n \to \infty$, $\left\| \xi_n \right\|_\infty \leq R$, $\overline{\xi_n}=0$, $\Phi_{r_n,m}(\xi_n) \in W^{1,2}(\Omega;\nu,\mu)$ and
\begin{equation}\label{eq: negazione}
\left\| \nabla{\Phi_{r_n,m}(\xi_n)} \right\|_{2;\mu} \leq \frac{1}{n} \left\| \Phi_{r_n,m}(\xi_n) \right\|_{2\sigma;\nu} \, .
\end{equation}
The fact that the sequence $\{ r_n \}$ cannot accumulate at some fixed $r$ and therefore has to go to infinity is a straightforward consequence of the method of proof of \cite[Lem. 5.8]{GMP}. Now, let us set
$$ a_n=\| \Phi_{r_n,m}(\xi_n) \|_{{2\sigma;\nu}^{\phantom{A}}} $$
and
$$ \Psi_n = \frac{\Phi_{r_n,m}(\xi_n)}{a_n} \, .  $$
From \eqref{eq: negazione} we get that
$$ \left\| \Psi_n \right\|_{2\sigma;\nu}=1 \, , \  \left\| \nabla{\Psi_n} \right\|_{2;\mu}\le\frac{1}{n} \, . $$
Applying the Poincar\'e-Sobolev inequality \eqref{eq: poinSob} to $\Psi_n$ we obtain
\begin{equation}\label{eq: poinSobApp}
\left\| \Psi_n - \overline{\Psi}_n \right\|_{2\sigma;\nu} \leq  C_{PS} \frac{1}{n} \, ;
\end{equation}
in particular, since $\left\| \Psi_n \right\|_{2\sigma;\nu}=1$, the sequence of real numbers $\left\{ \overline{\Psi}_n \right\}$ is bounded and therefore convergent to some constant $c_0$ (up to a subsequence which we do not relabel). Hence by \eqref{eq: poinSobApp} also $\{ \Psi_n \}$ converges in $L^{2\sigma}(\Omega;\nu)$ to such constant $c_0$: in particular $c_0 \neq 0$ since $\| \Psi_n \|_{2\sigma;\nu}=1$.
Consider the sequence of functions
$$ \mathcal{Z}_n = \frac{\xi_n}{a_n^{\frac{1}{r_n}}} \, . $$
First of all, we want to prove that $\{ \mathcal{Z}_n \}$ converges at least pointwise to a non-zero constant. Indeed, again up to subsequences, we know that $\{ \Psi_n \}$ converges pointwise to the non-zero constant $c_0$. Let us rewrite $\mathcal{Z}_n$ in the following way:
\begin{equation*}\label{eq: znRewrite}
\mathcal{Z}_n(\mathbf{x})= \left[ \frac{\left| \xi_n(\mathbf{x}) \right|^{r_n}}{\left| \Phi_{r_n,m}(\xi_n(\mathbf{x})) \right|}   \Psi_n(\mathbf{x}) \right]^{\frac{1}{r_n}} \, .
\end{equation*}
Thanks to estimate \eqref{eq: lowUppBound} and to the fact that $|\xi_n(\mathbf{x})| \leq R$, we have:
$$  \left[\frac{r_n}{C_2}\right]^{\frac{1}{r_n}} \left|\Psi_n(\mathbf{x}) \right|^{\frac{1}{r_n}}   \leq  |\mathcal{Z}_n(\mathbf{x})| \leq  \left[ \frac{r_n^{1+ 1 \vee \left[\frac{m-1}{2}\right] }}{C_1} \right]^{\frac{1}{r_n}} \left|\Psi_n(\mathbf{x}) \right|^{\frac{1}{r_n}} \, .  $$
Therefore by letting $n\to \infty$ we deduce that $\{ |\mathcal{Z}_n| \}$ converges pointwise to $1$, and so $\{ \mathcal{Z}_n\}$ converges pointwise to $1$ if $c_0 > 0 $ or to $-1$ if $c_0<0$, in any case to a non-zero constant. Finally we prove that $\{\mathcal{Z}_n\}$ also converges in $L^1(\Omega;\nu)$ to such non-zero constant. In order to do that, thanks to Egoroff's Theorem it is enough to show that $\int_E |\mathcal{Z}_n| \, \mathrm{d}\nu $ goes to zero as $n \rightarrow \infty$ and $|E| \rightarrow 0$. In fact,
$$ \int_E |\mathcal{Z}_n| \mathrm{d}\nu = \int_E { \frac{|\xi_n|}{a_n^{\frac{1}{r_n}}}   \mathrm{d}\nu } \leq \left[\frac{r_n^{1+ 1 \vee \left[\frac{m-1}{2}\right] }}{C_1} \right]^{\frac{1}{r_n}} \int_E \left|\Psi_n \right|^{\frac{1}{r_n}}  \mathrm{d}\nu  \leq |E|^{1-\frac{1}{r_n}}  \left[\frac{r_n^{1+ 1 \vee \left[\frac{m-1}{2}\right] }}{C_1} \right]^{\frac{1}{r_n}}  \left\| \Psi_n \right\|_{1;\nu}^{\frac{1}{r_n}} , $$
and the assertion follows since $\{\Psi_n\}$ converges in $L^{2\sigma}(\Omega;\nu)$, so in particular is bounded in $L^1(\Omega;\nu)$. Thus $\{\mathcal{Z}_n\}$ is a sequence of zero-mean functions which converge in $L^1(\Omega;\nu)$ to a non-zero constant, and this is clearly absurd. Therefore the constant $C_\ast$ in \eqref{eq: poinSobPhi} must exist.
\end{proof}
\end{lem}
A crucial point in Lemma \ref{lem: funct-inequality} lies in the fact that the constant $C_\ast$ does not depend on $r$, and this is necessary in order to prove the following theorem, which improves considerably the main result of \cite{KR2}, at least in all the cases dealt with therein for which a weighted Sobolev inequality holds. Indeed, in \cite{KR2} \it local \rm uniform convergence is shown for the specific weights considered in the one-dimensional setting. For such weights a Sobolev-type inequality may or may not hold, and hence one does not expect uniform convergence in general. Our contribution is the identification of the role of a weighted Sobolev inequality in this connection and the functional analytic setting in which we prove our results. 
\begin{thm}\label{teo: convUnif}
Let $\nu(\Omega)<\infty$ and let the Sobolev-type inequality \eqref{eq: soboIntroMedWei} hold true for some $\sigma>1$. Then any solution $u$ of \eqref{eq: pmeNeuWei} corresponding to an initial datum $u_0 \in L^{q_0}(\Omega;\nu)$, with $q_0 \in [1,\infty)$ and $\overline{u_0}=\overline{u}\neq 0$, converges \emph{uniformly} to its mean value $\overline{u}$.
\begin{proof}
We shall proceed through a Moser iterative technique. Indeed, for a given $\varrho>1$, after multiplying equation \eqref{eq: pmeNeuWei} by $(u-\overline{u})^{\varrho-1}$ and (formally) integrating by parts, straightforward computations lead to the following differential equation for the quantity $\| w(s) \|_{\varrho;\nu}^\varrho $, where we set $w$ to be the relative error $w(\cdot,s)=u(\cdot,s)/\overline{u}-1$:
\begin{equation}\label{eq: diffErr}
\frac{\mathrm{d}}{\mathrm{d}s} \left\| w(s) \right\|_{\varrho;\nu}^\varrho = -\varrho(\varrho-1)m|\overline{u}|^{m-1} \left\| \nabla{\Phi_{\frac{\varrho}{2},m}(w(s))} \right\|_{2;\mu}^2 \, .
\end{equation}
Since $u(\tau) \in L^\infty(\Omega)$ for any given $\tau>0$ (recall the regularizing effect from Theorem \ref{teo: soboTeoRegWei})  and $\| w(s) \|_\infty$ is non-increasing (consequence of \eqref{eq: diffErr} itself), of course there exists a suitable $R>1$ such that $|w(\cdot,s)| \leq R $ for all $s \geq \tau$; notice that the constant $C_\ast$ in Lemma \ref{lem: funct-inequality} will then depend on the solution $u$ only through $\|w(\tau)\|_\infty$. Also, by definition, $\overline{w}(s)=0$. Therefore we can apply to the right hand side of \eqref{eq: diffErr} the functional inequality \eqref{eq: poinSobPhi} (with $r=\varrho/2$), which entails
\begin{equation}\label{eq: diffErr1}
\frac{\mathrm{d}}{\mathrm{d}s} \left\| w(s) \right\|_{\varrho;\nu}^\varrho \leq -\frac{ \varrho(\varrho-1)m|\overline{u}|^{m-1}}{C_\ast^2} \left\| \Phi_{\frac{\varrho}{2},m}(w(s)) \right\|_{2\sigma;\nu}^2 \, .
\end{equation}
Thanks to the estimate \eqref{eq: lowUppBound} (again with $r=\varrho/2$), from \eqref{eq: diffErr1} we get
\begin{equation}\label{eq: diffErr2}
\frac{\mathrm{d}}{\mathrm{d}s} \left\| w(s) \right\|_{\varrho;\nu}^\varrho \leq -\frac{ \varrho(\varrho-1)m|\overline{u}|^{m-1}}{C_\ast^2} \frac{C_1^2}{\left( \frac{\varrho}{2} \right)^{2+2\vee(m-1)}}  \left\| w(s)^{\frac{\varrho}{2}} \right\|_{2\sigma;\nu}^2 \leq - \frac{Q}{\varrho^{2\vee(m-1)}} \left\| w(s) \right\|_{ \sigma \varrho ;\nu }^{\varrho} \, ,
\end{equation}
where $Q>0$ is a suitable constant depending on $m$, $C_1$, $C_\ast$, $|\overline{u}|$,  but independent of $\varrho$. Now, for a given $t>\tau$, let us set $t_n=\tau+(1-2^{-n})(t-\tau) $ and $p_{n+1}=\sigma \,  p_n$. Clearly, $p_n=\sigma^n p_0$ (let $p_0$ be any real number belonging to $[1,\infty)$), so that we have $t_n \to t$ and $p_n \to \infty$ as $n \to \infty$. Integrating \eqref{eq: diffErr2} between $t_n$ and $t_{n+1}$ with $\varrho=p_n$ and exploiting the non-expansivity of $\| w(s) \|_{p_{n+1};\nu}$ we obtain:
$$ -\left\| w(t_n) \right\|_{p_n;\nu}^{p_n} \leq - \frac{Q}{p_n^{2\vee(m-1)}} (t_{n+1}-t_n) \left\| w(t_{n+1}) \right\|_{p_{n+1};\nu}^{p_n}  \, ,$$
which reads
\begin{equation}\label{eq: disRic}
\left\| w(t_{n+1}) \right\|_{p_{n+1};\nu} \leq \left(\frac{ p_n^{2\vee(m-1)} \, 2^{n+1}}{Q}  \right)^{\frac{1}{p_n}} \, t^{-\frac{1}{p_n}} \, \left\| w(t_n) \right\|_{p_n;\nu} \, .
\end{equation}
Since $p_n=\sigma^n p_0$, for a suitable constant $K=K(p_0,m,\sigma,Q)$ we can rewrite \eqref{eq: disRic} as
\begin{equation*}\label{eq: disRic2}
\left\| w(t_{n+1}) \right\|_{p_{n+1};\nu} \leq K^{\frac{n+1}{p_n}} \, t^{-\frac{1}{p_n}} \, \left\| w(t_n) \right\|_{p_n;\nu} \, ,
\end{equation*}
that is
\begin{equation*}\label{eq: disRic3}
\left\| w(t_{n+1}) \right\|_{p_{n+1};\nu} \leq K^{\sum_{k=0}^n { \frac{k+1}{p_k} } } \, t^{-\sum_{k=0}^n \frac{1}{p_k} } \, \left\| w(\tau) \right\|_{p_0;\nu} \, .
\end{equation*}
Hence,
\begin{equation}\label{eq: disFinale}
\left\| w(t) \right\|_\infty = \lim_{n\to \infty} \left\| w(t) \right\|_{p_{n+1};\nu} \leq \liminf_{n\to \infty} \left\| w(t_{n+1}) \right\|_{p_{n+1};\nu} \leq K^\prime \, t^{- \frac{\sigma}{p_0(\sigma-1)} } \, \left\| w(\tau) \right\|_{p_0;\nu} \, ,
\end{equation}
being $K^\prime$ another suitable constant depending $K$, $p_0$, $\sigma$ and therefore on $p_0$, $m$, $\Omega$, $\nu$, $\mu$, $\overline u$, $\|w(\tau)\|_\infty$. The assertion then follows by letting $t\to \infty$ in \eqref{eq: disFinale}.
\end{proof}
\end{thm}
The above result can be refined by showing the following analogue of Theorem \ref{eq: teoAsiMNon}.
\begin{thm}\label{th: teoAsiMNonPesi}
Let $\nu(\Omega)<\infty$ and let the Sobolev-type inequality \eqref{eq: soboIntroMedWei} hold true for some $\sigma>1$. For any solution $u$ of \eqref{eq: pmeNeuWei} corresponding to an initial datum $u_0 \in L^1(\Omega;\nu)$ with $\overline{u_0}=\overline{u}\neq 0$ there exists a constant $G=G(m, \Omega, \nu, \mu, \overline{u}, \| u_0\|_{1;\nu})>0$ such that the following estimate holds:
 \begin{equation*}\label{eq: StimaExpFinWei}
\left\| u(t)-\overline{u} \right\|_\infty \leq G \, e^{-\frac{m}{C_P^2}|\overline{u}|^{m-1} \, t   } \ \ \ \forall t \geq 1  \, ,
\end{equation*}
being $C_P$ the smallest constant such that \eqref{eq: poinIntroMedWei} holds.
\begin{proof}
One can proceed along the lines of the proof of Theorem \ref{eq: teoAsiMNon}. In order to get the expected exponential decay of the quantity $\| u(t)-\overline{u} \|_{2;\nu}$, one first uses estimate \eqref{eq: disFinale} instead of \eqref{eq: stimaConvUnif2}. Indeed, note that the only relevant point of \eqref{eq: stimaConvUnif2} is that its right hand side goes to zero and it is \emph{integrable} as $t\to \infty$: both of these properties can be achieved by the right hand side of \eqref{eq: disFinale} upon choosing $p_0$ sufficiently next to $1$ (or even $p_0=1$). Finally, the decay rate of the $L^\infty$ norm is the same as for the $L^2$ norm thanks again to \eqref{eq: disFinale} evaluated (for example) between $t$ and $\tau=t-1/2$ with $p_0=2$. The dependence of the multiplicative constant $G$ on the stated quantities follows from the constant $K^\prime$ in \eqref{eq: disFinale} and from the regularity estimate \eqref{teo: soboTeoRegWei}.
\end{proof}
\end{thm}
We remark that, obviously, Theorem \ref{th: teoAsiMNonPesi} also applies to the non-weighted case. In fact the conclusion is the same as the one of Theorem \ref{eq: teoAsiMNon} in terms of time decay rate. Nonetheless, in Theorem \ref{eq: teoAsiMNon} the multiplicative constant $G$ depends in a nontrivial way on the solution itself, while the proof of Theorem \ref{th: teoAsiMNonPesi} shows that, basically, it depends on the solution only in terms of the initial datum (though through the constant $C_\ast$ of Lemma \ref{lem: funct-inequality}, which is unknown).

We conclude giving an immediate corollary of Theorem \ref{teo: convUnif} concerning the evolution of the support of (non-zero mean) solutions to \eqref{eq: pmeNeuWei}. This topic has been widely investigated in the literature, see e.g. \cite{KK, P, GHP, KRT, T}. 
\begin{cor}\label{cor: support}
Let $\nu(\Omega)<\infty$ and let the Sobolev-type inequality \eqref{eq: soboIntroMedWei} hold true for some $\sigma>1$. Then the support of any solution $u(\cdot,t)$ of \eqref{eq: pmeNeuWei} corresponding to a \emph{compactly supported} initial datum $u_0$ having non-zero mean becomes the whole $\Omega$ for all $t$ great enough.
\end{cor}
We recall that in \cite{GMP}, where the weighted equation \eqref{eq: pmeNeuWei} was analysed assuming only the validity of the Poincar\'e inequality \eqref{eq: poinIntroMedWei}, the authors provided weights for which the corresponding solutions to \eqref{eq: pmeNeuWei} do not converge uniformly to their mean value (but they do in $L^\varrho(\Omega;\nu)$ for all $\varrho \in [1,\infty)$ as a consequence of \eqref{eq: poinIntroMedWei} itself). This should indicate that \emph{uniform} convergence is strongly linked with the validity of \eqref{eq: soboIntroMedWei}.
\subsection{Examples of weighted Sobolev inequalities }\label{sec: exaSobo}
In the following we list some basic examples of domains $\Omega$ and couples of weights $(\rho_\nu,\rho_\mu)$ for which the Sobolev-type inequality \eqref{eq: soboIntroMedWei} holds for appropriate values of the parameter $\sigma>1$, given below.
\begin{itemize}
\item[$\bullet$] \emph{Intervals}:
\begin{itemize}
\vskip 1mm
\item[$\circ$] $(x^{\alpha},x^\beta)$ on $(0,b)$ ($b>0$): $\beta \leq 1$, $\alpha > -1 $ and all $\sigma > 1 $ OR $\beta > 1$, $\alpha > \beta-2 $ and $\sigma\in \left(1,\frac{\alpha+1}{\beta-1}\right] $;
\item[$\circ$] $(x^{\alpha},x^\beta)$ on $(a,+\infty)$ ($a>0$): $\beta \geq 1$, $\alpha < -1 $ and all $\sigma > 1 $ OR $\beta < 1$, $\alpha < \beta-2 $ and $\sigma\in \left(1,\frac{\alpha+1}{\beta-1} \right] $;

\vskip 1mm
\item[$\circ$] $\left( \frac1{x}|\log{x}|^{\alpha}, {x}|\log{x}|^{\beta} \right)$ on $(0,c)$ ($c\in(0,1)$): $\beta \geq 1$, $\alpha < -1 $ and all $\sigma > 1 $ OR $\beta < 1$, $\alpha < \beta-2 $ and $\sigma\in \left(1, \frac{\alpha+1}{\beta-1} \right] $;

\vskip 1mm
\item[$\circ$] $(e^{\alpha |x|},e^{\beta |x|})$ on $\mathbb{R}$: $\beta \geq 0$, $\alpha<0$ and all $\sigma>1$ OR $ \beta <0 $, $\alpha < \beta$ and $\sigma\in \left( 1, \frac{\alpha}{\beta} \right] $.
\end{itemize}

\vskip 1mm
\item[$\bullet$] \emph{Bounded Lipschitz domains} ($N>1$):
\begin{itemize}
\vskip 1mm
\item[$\circ$] $(\delta^{\alpha},\delta^\beta)$: $\beta \leq 1$, $\alpha > -1 $ and $\sigma \in \left(1, \min\left( \frac{N}{N-2}, \frac{\alpha+N}{N-1} \right) \right] $ OR $\beta > 1 $, $\alpha>\beta-2$ and $\sigma \in \left(1, \min\left( \frac{N}{N-2}, \frac{\alpha+N}{\beta+N-2} \right) \right] $  ($\delta$ denotes the distance function from $\partial\Omega$).
\end{itemize}
\vskip 1mm
\item[$\bullet$] \emph{The Euclidean space $\mathbb{R}^N$} ($N>1$):
\begin{itemize}
\vskip 1mm
\item[$\circ$] $((1+|\mathbf{x}|)^{\alpha},(1+|\mathbf{x}|)^{\beta})$: $ \beta \geq 2-N$, $\alpha<-N $ and $\sigma \in \left(1,\frac{N}{N-2} \right] $ OR $ \beta < 2-N $, $\alpha < \beta-2 $ and $\sigma \in \left(1, \min\left(\frac{N}{N-2} , \frac{\alpha+N}{\beta+N-2} \right) \right] $;
\item[$\circ$] $( e^{\alpha|\mathbf{x}|} , e^{\beta|\mathbf{x}|} )$: $\beta \geq 0$, $\alpha < 0$ and $\sigma \in \left(1,\frac{N}{N-2} \right]$ OR $\beta < 0$, $\alpha<\beta$ and $\sigma \in \left( 1, \min\left(\frac{N}{N-2},\frac{\alpha}{\beta} \right) \right] $.
\end{itemize}
\end{itemize}
The above examples can be obtained by applying Theorem 1.4 of \cite{CW99} (which gives necessary and sufficient conditions on the weights so that \eqref{eq: soboIntroMedWei} holds) in the one-dimensional case and by exploiting the results of \cite{KO90} in the $N$-dimensional case, see in particular Chapter 19 for bounded Lipschitz domains and Chapters 20, 21 for unbounded domains (the validity of \eqref{eq: soboIntroMedWei} here is a consequence of the compact embeddings discussed therein).
\end{section}

\bibliographystyle{plainnat}


\end{document}